\documentstyle[12pt]{article}

 \title{ On Permanental Processes}
 \date{}
\newtheorem{Theorem}{Theorem}[part]
\newtheorem{Definition}[Theorem]{Definition} 
\newtheorem{Proposition}[Theorem]{Proposition}

\newtheorem{Lemma}[Theorem]{Lemma} 
\newtheorem{Corollary}[Theorem]{Corollary}

\makeatletter \@addtoreset{equation}{section}

\@addtoreset{Theorem}{section}

\font\tenmath=msbm10
\font\sevenmath=msbm7
\font\fivemath=msbm5
\newfam\mathfam \textfont\mathfam=\tenmath
\scriptfont\mathfam=\sevenmath
\scriptscriptfont\mathfam=\fivemath
\def\math{\fam\mathfam}
\def \b{\noindent}
\def \e{I\!\!E}
\def \p{I\!\!P}

\def \={{\buildrel {\rm (law)} \over =}}
\def \N{{\math N}}
\def \R{{\math R}}

\def \1{1\!\!1}
\def \l{\ell}

\begin{document}

\maketitle

\centerline{\bf Nathalie Eisenbaum and Haya Kaspi}

\bigskip

\bigskip
\bigskip

\b {\bf Abstract} : Permanental processes can be viewed as a generalization of squared centered Gaussian processes. We develop in this paper two main directions.  The first one analyses the connections of these processes with the local times of general Markov processes. The second  deals with Bosonian point processes and the Bose-Einstein condensation. The obtained results in both directions are related and based on the notion of infinite divisibility.

\bigskip

\bigskip

\b {\bf Keywords} : Gaussian process, local times, infinite divisibility, permanental process, random point process, Bose-Einstein condensation. 

\bigskip

\bigskip

\b  {\bf AMS 2000 subject classification :} 60J25, 60J55, 60G15, 60G55

\bigskip

\vfill
\eject

 \section{Introduction}
 
 \bigskip
 \b Permanental processes can be viewed as a generalization of the squared centered Gaussian processes. Their Laplace transform is given by the  power $(-{1\over \alpha})$ of a determinant  ($\alpha > 0$) involving a kernel,  squared Gaussian processes corresponding to the case of a symmetric kernel and $\alpha = 2$. The value of $\alpha$ is called the index of the permanental process. The problem of the existence of such processes has been solved by Vere-Jones \cite{VJ}.   This paper develops mostly two subjects related to these processes. The first one analyses  connections between permanental processes and the local times of a general Markov process. The second one deals with permanental random point processes  also called Bosonian random point processes.

\b  The first subject is based on the natural mergence of the permanental processes in the study of the local times of Markov processes.  In the case of a symmetric Markov process, this presence has allowed the writing of so-called 
``isomorphism theorems" connecting directly the law of the local times to the law of a squared Gaussian process. The most famous one is the identity of Dynkin \cite{D1}.   Marcus and Rosen's book \cite{MR1} makes obvious  the interest of these theorems. In the general case (non necessarily symmetric), a permanental process is going to replace this squared Gaussian process and identities similar to Dynkin's isomorphism theorem can then be written.  We establish here two identities : one  for the total accumulated local time of a transient Markov process, and a second one for recurrent Markov processes stopped at inverse local times (extending an identity of \cite{gang}).

\b The problem  then is to be able to use these identities.  As an example of use, we show here that if the local time process of a transient Markov process is continuous then  the  associated permanental  process is continuous.  To do so, 
we use the fact  that the permamental processes  associated to Markov processes, are always infinitely divisible. 
 
 \b In previous works, we have shown \cite{E1},\cite{EK} that the property of infinite divisibility characterizes the squared Gaussian processes associated to Markov processes. Here we extend this characterization to the non-symmetric case: namely we show that a  permanental process is infinitely divisible iff it is associated to a Markov process. Moreover, its L\'evy measure is given. 
 
 \bigskip

 \b The definition of permanental  (or alpha-permanental)  random point processes is due to Shirai and Takahashi \cite{ST}. Indeed they have established the existence of random point processes such that their  Laplace transform is equal to the power $(-{1 \over \alpha})$ of a Fredholm determinant. When $\alpha = -1$ one obtains a determinantal  (or a fermion) random process, when $\alpha = 1$ it is a boson point process.  When $\alpha$ is positive one can call them permanental random processes because their densities and their correlation functions are equal to permanents.
  Besides when $\alpha$ is positive, there is a  connection between the corresponding  permanental random point process and the real permanental process with index $\alpha$ under the condition that this last exists. Indeed, in that case the permanental random point process is a Cox process driven by a permanental process with index $\alpha$.

\b The most known illustration of these processes, is the law of the configurations of an ideal gaz of Bosonian particles in standard conditions (see for example Macchi \cite{Ma}). Recently Tamura and It\^o \cite{TI} have obtained  the law of the configurations of the particles of an ideal Bosonian gaz  containing 
particles in a Bose-Einstein condensation state.  We analyze here their result and
show that the obtained law provides an illustration of some kind of "super"- Isomorphism Theorem existing above the usual isomorphism theorem. 

\bigskip

\b The paper is organized as follows.  Section 2 treats about general properties of permanental processes such as existence, conditioning, absolute continuity.  Section 3 specializes in permanental processes associated to Markov processes. Section 4 provides a characterization of the infinitely divisible permanental processes.  Section 5 deals with permanental random point processes and analyzes the result of Tamura and Ito concerning the Bose-Einstein condensation phenomena. Section 6 provides a factorization result on positive infinitely divisible processes. Section 7 contains the proof of the results exposed in Sections 2, 3 and 4.  We end the paper by a translation of Shirai and Takahashi's  conjecture on $\alpha$-permanents  in Section 8. 

\bigskip

\b We mention that a recently posted on ArXiv paper by Yves Le Jan \cite{LJ}, contains a result similar to the  isomorphism theorem established here.

 \section{Existence, conditioning and behavior} \label{S1}

 \b By permanental processes, we mean processes such that their Laplace transform is given by a negative power of a determinant.  More precisely 
 \begin{Definition}\label{001}: A real-valued positive process $(\psi_x , x\in E)$ is a permanental process if its finite-dimensional Laplace transforms satisfy
 for every 
$(\alpha_1, \alpha_2, ...,\alpha_n)$ in $\R^n_+$ and every $(x_1,x_2, ...,x_n)$ in $E^n$,
\begin{equation}
\e[\hbox{exp}\{-{1\over 2} \sum_{i = 1}^n \alpha_i \psi_{x_i}\}] = | I + \alpha G |^{-1/\beta}
\end{equation}
where $I$ is the $n\times n$-identity matrix, $\alpha$ is the diagonal matrix $\hbox{diag}(\alpha_i)_{1 \leq i \leq n}$ and 
$G = (g(x_i,x_j))_{1\leq i,j \leq n}$ and $\beta$ is a fixed positive number.

\b Such a process $(\psi_x , x\in E)$ is called permanental process with kernel 

\b $(G(x,y), x,y \in E)$ and index $\beta$.
\end{Definition}
Vere-Jones has established in \cite{V} the above necessary and sufficient conditions (1) and (2) on matrices $G = (g(x_i,x_j))_{1\leq i,j \leq n}$ for the existence of a corresponding permanental vector $(\psi_{x_i}, 1 \leq i \leq n)$. They are based on the following definition.

\begin{Definition} : For any $n\times n$ matrix M : $$det_{\beta} M = \sum_{\sigma \in {\cal S}_n} \beta^{n - \nu(\sigma)}\Pi_{i = 1}^n M_{i,\sigma(i)},$$ where ${\cal S}_n$ is the symmetric group of order $n$  and $\nu(\sigma)$ is the number of cycles of 
$\sigma$. 
For every multi-index $k = (k_1,k_2,...,k_n)$ $M(k)$ denotes the derived  $|k|\times |k|$-matrix (where $|k| = k_1+k_2+...+k_n$) obtained from $M$ by selecting the first row and column $k_1$ times, the second $k_2$ times,..., the $n^{th}$ $k_n$ times. For $\beta > 0$, a matrix $M$ is said to be  {\bf $\beta$-positive definite} if for all possible derived matrices $M(k)$, $det_{\beta} (M(k)) \geq 0$.
\end{Definition} 
A permanental vector $(\psi_{x_i}, 1 \leq i \leq n)$ corresponding to $G = (g(x_i,x_j))_{1\leq i,j \leq n}$ and index $\beta$  exists if and only if:

{\bf (1)} $|I + r G| > 0$ for every $r>0$.

{\bf (2)} For every $r>0$, set $Q_r = G(I + rG)^{-1}$, then $Q_r(k)$ is $\beta$- positive definite for every $k$ in $\N^n$.

\bigskip

\b Here is our first result concerning these permanental processes. 
The proof is provided in Section \ref{S6}.

\begin{Proposition}\label{P4}: Let $(G(x,y), (x,y) \in E\times E)$ be a real function on $E\times E$ such that there exists $a$ in $E$ with $G(x,a)= G(a,x) = 0$ for every x in $E$. 
For a fixed  $\delta > 0$ assume that there exists a permanental process $(\psi_x, x \in E)$ with a kernel $(G(x,y) + \delta, (x,y) \in E\times E)$ and index $\beta > 0$.  We have then
\begin{equation}
 \label{01}
 \e[\hbox{exp}\{-{t\over 2}  \psi_a\}] =   (1 + \delta t)^{-1/\beta} 
\end{equation}
and for every $n$,  every 
$(\alpha_1, \alpha_2, ...,\alpha_n)$ in $\R^n_+$, every $(x_1,x_2, ...,x_n)$ in $E^n$,  every $r\geq 0$
\begin{equation}
\label{02}
 \e[\hbox{exp}\{-{1\over 2} \sum_{i = 1}^n 
 \alpha_i \psi_{x_i}\} | \psi_a = r] = | I + \alpha G |^{-1/\beta} \hbox{exp}\{- {1\over 2}r 1^t(I + \alpha G)^{-1}\alpha 1 \} 
\end{equation}
where  $I$ is the $n\times n$-identity matrix, $\alpha$ is the diagonal matrix $\hbox{diag}(\alpha_i)_{1 \leq i \leq n}$ and 
$G = (g(x_i,x_j))_{1\leq i,j \leq n}$. 
\end{Proposition}
The existence of $\psi$ is equivalent to Vere-Jones conditions (1) and (2) for $(G + \delta)$, for every $n$ and every $x_1,x_2,...x_n$ in $E$.   As a consequence of Proposition \ref{P4}, we obtain, under the same assumptions, the following result which can not be easily seen using just (1) and (2).
\smallskip

\b {\sl  If there is a  $\delta_0 > 0$  such that there exists a permanental process $\psi$ with kernel $(G + \delta_0)$  and index $\beta> 0$, then for every $\delta \geq 0$ there exists a permanental process with kernel $(G + \delta)$ and index $\beta$. In particular, denoting by $\phi$ a  permanental process with kernel $G$ and index $\beta$, we obtain 
$$(\psi_x , x\in E \ | \ \psi_a = 0) \ \= \ (\phi_x , x\in E).$$}

\b In the case when $G$ is positive definite and $\beta = 2$, permanental processes with kernel $G + \delta$ and index $2$ exists for every $\delta \geq 0$ (see the remark below).  
Section \ref{S2} deals with a class of kernels $G$ for which permanental processes with kernel $(G + \delta)$  and index $\beta$ exist for every $\delta \geq  0$ and $\beta \geq 0$.  
\smallskip

\b Proposition \ref{P4} implies the following property for $\psi$
\begin{eqnarray*}
(\psi_x, x\in E | \psi_a=r) \ &+& \ (\tilde{\psi}_x, x\in E| \tilde{\psi}_a=r') \\
& \= &(\psi_x, x\in E | \psi_a=t) \ + \ (\tilde{\psi}_x, x\in E| \tilde{\psi}_a=t') 
\end{eqnarray*}
for $\tilde{\psi}$ an independent copy of $\psi$, and $r$,$r'$,$t$, $t'$ nonnegative numbers such that $r+r' = t+ t'$.
This property is well-known when $(G(x,y), (x,y) \in E\times E)$ is symmetric and $\beta = 2$. Indeed, in that case
 $(\phi_x, x\in E)$ is a squared centered Gaussian process. More precisely there exists a centered Gaussian process $(\eta_x , x \in E)$   with covariance $(G(x,y), x,y \in E)$ such that 
$(\phi_x, x\in E) = (\eta_x^2 , x \in E)$. One can always add an element $a$ by setting $G(a,x) = G(x,a) = 0$, $(\eta_x , x \in E \cup \{a\})$ remains a centered Gaussian process.
We have then: $\psi \= (\eta + N)^2$ where $N$ is a centered Gaussian variable with a variance equal to $\delta$, independent of $\eta$. This gives
$$(\psi_x, x \in E | \psi_a = r^2)\  \= \ ((\eta_x  + N)^2, x \in E  \ | N = r  ) \ \= \  ((\eta_x + r)^2 , x \in E)$$
 and for $\tilde{\eta}$  an independent copy of $\eta$, we have: 
\begin{eqnarray*}
 (\eta  + a)^2 + (\tilde{\eta}  + b)^2 \= (\eta  + c)^2 + (\tilde{\eta}  + d)^2.
\end{eqnarray*}
for any $a$,$b$,$c$ and $d$ such that $a^2 + b^2 = c^2 + d^2$.

\bigskip

\b Although  the shape of the Laplace transforms of a permanental process  is close to the one of a squared Gaussian process, there is no result in the literature  on the pathwise behavior of these processes. The following proposition  connects some permanental processes with index $2$ to  squared Gaussian processes and should give the key to understand  their pathwise behavior.

\begin{Proposition}\label{P5}: Let $(\psi_x , x \in E)$ be a permanental process with kernel $(G(x,y), (x,y) \in E\times E)$ and index $2$. Assume that  $({1\over 2}(G(x,y) + G(y,x)), (x,y) \in E\times E)$ is positive definite.  Let $(\eta_x, x\in E)$ be a centered Gaussian process with this covariance. Let $\tilde{\eta}$ be an independent copy of $\eta$ and $\tilde{\psi}$ an independent copy of $\psi$. 
For  $x_1,x_2,...,x_n$ in $E$,  we set $G = (G(x_k,x_j))_{1 \leq k,j \leq n}$, and  $\Lambda = (\eta_{x_k} + i \tilde{\eta}_{x_k}, 1 \leq k \leq n)$. 
We have the following relation for any functional $F$ on $\R^n$
\begin{equation}
\label{05}
 \e[ F(\psi_{x_k}^2 + \tilde{\psi}_{x_k}^2 , 1 \leq k \leq n)] = {\e[  \hbox{exp}\{{1\over 2}<A\Lambda,\Lambda>\} F(\eta^2_{x_k} + \tilde{\eta}^2_{x_k}, 1 \leq k \leq n)] \over \e[  \hbox{exp}\{{1\over 2}<A\Lambda,\Lambda>\} ]}
 \end{equation}
where $A_{ij} =   ({1\over 2}(G + G^t))^{-1} - G^{-1})_{ij}$.

\end{Proposition}
Of course if the kernel $G$ is symmetric, $\eta^2$ coincides with $\psi$. 

\bigskip
 
 \section{Permanental processes associated to Markov processes.}\label{S2}

\bigskip

\b We work with a transient Markov process with a state space $E$. Denote by $(L^x_t, x\in E, t\geq 0)$ its local time process and by $(g(x,y), (x,y)\in E\times E)$ its Green function. It satisfies : $g(x,y) = \e_x(L^y_{\infty})$. Let $a$ be an element of $E$ such that $g(a,a) > 0$. We define the probability $\tilde{\p}_{a}$  as follows 
$$\tilde{\p}_{a}|_{{\cal F}_t} = {g(X_t,a)\over g(a,a)}\p_a|_{{\cal F}_t}$$ 
where ${\cal F}_t$ denotes the field generated by $(X_s, 0\leq s\leq t)$ and $\p_a$
the probability under which $X$ starts at $a$. 

\b Under $\tilde{\p}_{a}$, the process $X$ starts at $a$ and is killed at its last visit to $a$. Expectation with respect to $\tilde{\p}_{a}$ is denoted by $\tilde{\e}_{a}$.

\bigskip

\begin{Theorem}\label{T1} For every $\beta > 0$, 
 there exists a positive  process $(\psi_x , x\in E)$ 
such that for every 
$(\alpha_1, \alpha_2, ...,\alpha_n)$ in $\R^n_+$ and every $(x_1,x_2, ...,x_n)$ in $E^n$,
\begin{equation}
\e[\hbox{exp}\{-{1\over 2} \sum_{i = 1}^n \alpha_i \psi_{x_i}\}] = | I + \alpha G |^{-1/\beta}
\end{equation}
where $I$ is the $n\times n$-identity matrix, $\alpha$ is the diagonal matrix $\hbox{diag}(\alpha_i)_{1 \leq i \leq n}$ and 
$G = (g(x_i,x_j))_{1\leq i,j \leq n}$. 
\end{Theorem}

\bigskip

\b In the case $\beta = 2$, note that for every fixed $x\in E$,  $\psi(x)$ has the law of squared centered Gaussian variable with a variance equal to $g(x,x)$. If moreover the Green function is symmetric, $\psi$ is the square of a centered Gaussian process with a covariance equal to $(g(x,y), x,y\in E)$. This has been already noted and exploited by many authors (Dynkin \cite{D1},\cite{D2}, Marcus and Rosen \cite{MR}, Eisenbaum \cite{E1}, Eisenbaum et al \cite{gang}, ...). This Gaussian process has been called the "Gaussian process associated" to $X$. 

\smallskip

\begin{Definition} In the general case, we call the process $\psi$, the {\bf permanental process with index $\beta$ associated} to $X$. 
\end{Definition}

\bigskip
\b We will see in Section \ref{S3} that even when the Green function $g$ is not symmetric, it might happend that the associated permanental process $\psi$ with index $2$ is a squared of a centered Gaussian process.

\bigskip

\b Let $(\psi_x , x\in E)$ be the permanental process with index $2$ associated to the Markov process $X$, defined on a probability space unrelated with  $X$. On this probability space, the expectation will be denoted by $\langle \  ; \ \rangle$. The following theorem provides a connection 
between the law of $\psi$ and the law of $(L^x_{\infty}, x \in E)$.
\bigskip

\begin{Theorem}\label{T2} 
: For every $a\in E$ such that $g(a,a) > 0$, for every functional $F$  on the space of measurable functions from $E$ to $\R$, we have 
\begin{equation}
\tilde{\p}_{a}\langle F(  L^x_{\infty}  + {1\over 2} \psi_x ; x\in E) \rangle = \langle  {\psi_a\over g(a,a)} F({1\over 2} \psi_x ; x\in E) \rangle
\end{equation}
\end{Theorem}
The existence  of an associated permanental process with index $\beta$ for every $\beta > 0$, provides immediately the property of infinite divisibility of these processes. But we also have their L\'evy measures.

\begin{Corollary}\label{C1}:  The process $\psi$ is an infinitely divisible process with a L\'evy measure $\nu$ characterized by 
the following marginals
\begin{eqnarray*}
\nu_{(\psi(a)/2,\psi(x_2)/2,...\psi(x_n)/2)}\!\!\!\!\!\!&&\!\!\!\!\!\!(dy_1,dy_2,...,dy_n) \\
&=& {g(a,a)\over 2y_1} \tilde{\p}_{a}(L^a_{\infty}\in dy_1,
L^{x_i}_{\infty}\in dy_i, 2\leq i\leq n) .
\end{eqnarray*}
Equivalently the L\'evy measure of $\psi/2$ is equal to the law of $(L^x_{\infty}, x\in E)$ under 
${g(a,a)\over 2L^a _{\infty}} \tilde{\p}_{a}$ for every $a$ in the state space.
\end{Corollary}
The above corollary provides interesting connections between the path properties of the process $(\psi_x, x\in E)$ and the path properties of the local time process. The following theorem is an immediat application of Corollary \ref{C1}. We assume that $(E,{d})$ is a locally compact metric space. We say that the local time  $L$ is ${d}$-continuous when it satisfies
a.s. for every $x$ in $E$ and every $t>0$  $$\lim_{\begin{array}{l}|s-t| \rightarrow 0 \\ {d}(x,y) \rightarrow 0\end{array}} L^y_s = L^x_t.$$

\begin{Theorem}\label{T21}: If the local time process $L$ is ${d}$-continuous , then $\psi$ is ${d}$-continuous.
\end{Theorem}

\b In the case when the Markov process $X$ is symmetric, we obtain hence that the associated Gaussian process is continuous. A result which has already been established by Marcus and Rosen \cite{MR}. Their proof is based on the isomorphism Theorem of Dynkin.  But we would like to emphasize  the shortcut provided by the property of infinite divisibility of $\psi$.

\b We end this section by a version in the  non symmetric case of  a theorem established in \cite{gang}. 
Assume that $X$ is a recurrent Markov process with a state space $E$. For $a \in E$, define $T_a = \inf\{t\leq 0 : X_t = a\}$ and  $\tau_r = \inf\{t\geq 0 : L^a_t > r\}$.  Let $S_{\theta}$ be an exponential time with parameter $\theta$, independent of $X$. Then  $X$ killed at $T_a$ and $X$ killed at $\tau_{S_{\theta}}$ are both transient Markov processes. We denote by $\phi$ and $\psi$ their respective associated permanental  processes with index $2$.
We have  the following identity for the process $(L^x_{\tau_r}, x\in E)$.

\bigskip

\begin{Corollary}\label{C2}  :  Let $X$ be a recurrent Markov process. For $a\in E$ and  
every functional $F$  on measurable function from $E$ to $\R$, we have 
\begin{equation}
{\p}_{a}\langle F(  L^x_{\tau_r}  + {1\over 2} \phi_x ; x\in E) \rangle = \langle   F({1\over 2} \psi_x ; x\in E) | \psi_a = r \rangle.
\end{equation}
Besides, we have : $( \psi_x ; x\in E) | \psi_a = 0 ) \= (\phi_x ; x\in E)$.

\end{Corollary}

 \section{ Characterization of the infinitely divisible permanental processes}\label{S3}
\bigskip

\b Similarly to what has been done in  the symmetric case (see \cite{EK}), one might ask whether this property of infinite divisibility characterizes the associated permanental processes.  The answer is affirmative :  a permanental process is infinitely divisible if and only if it is associated to a Markov process.  In particular, a squared Gaussian process is infinitely divisible if and only if
it is  a permanental process associated to a Markov process (this does not imply necessarely that the Gaussian process itself is associated to a Markov process).  If a permanental process with index $\beta > 0$ is infinitely divisible then the permanental process with the same kernel and index $2$ is infinitely divisible too, hence from now on in this section we will take $\beta = 2$. 

\b To present the proof of that answer we first establish the following criterion which represents an extension of Bapat's criterion without assumption of symmetry.
 \bigskip

\begin{Definition} :  A $n\times n$-matrix $A$ is an $M$-matrix if and only if

(i) $A_{ij} \leq 0$ for $i\not= j$

(ii) $A$ is non-singular and $A^{-1} \geq 0$ (i.e.  $A^{-1}_{ij} \geq 0$ for every $i,j$).
\end{Definition}
 
\bigskip
\begin{Lemma}\label{L1} : Let $(G_{i,j} ,1\leq i,j\leq n)$ be a real non-singular $n\times n$-matrix.  There exists a positive infinitely divisible random vector $(\psi_1, \psi_2, ...\psi_n)$ such that for every $(\alpha_1, \alpha_2, ..., \alpha_n) \in \R^n_+$,
\begin{equation}
\label{03}
\e[\hbox{exp}\{-{1\over 2} \sum_{i = 1}^n \alpha_i \psi_{i}\}] = | I + \alpha G |^{-1/2}
\end{equation}
if and only if, 
there exists a signature matrix $S$ such that $S G^{-1} S$ is an $M$-matrix.
\end{Lemma}
\bigskip

\b Note the following   consequence of the previous lemma :  the real eigenvalues of a matrix $G$ satisfying (\ref{03}) must be positive.

\begin{Theorem}\label{T4}: Let $(G_{i,j} ,1\leq i,j\leq n)$ be a real non-singular $n\times n$-matrix.  There exists a positive infinitely divisible random vector $(\psi_1, \psi_2, ...\psi_n)$ such that (\ref{03}) is satisfied, if and only if 
\begin{equation}\label{0}
G(i,j) = d(i)g(i,j)d(j)
\end{equation}
for every $(i,j)$,  where $d$ is a  function  on $\{1,2,...n\}$ and $g$ the Green function of a Markov process.
\end{Theorem}
\bigskip
\b {\bf Remark \ref{T4}.1} :  As it has been noticed in \cite{EK}, the property (\ref{0}) is equivalent to the following property
\begin{equation}
\label{15}
G(i,j) = d^{-1}(i)g(i,j)d(j)
\end{equation}
where $d$ is a  function and $g$ the Green function of a Markov process.

\b But then note that: $|I + \alpha G| = |I + \alpha g|$.  This means that the vector $\psi$ is a permanental vector associated to a Markov process.  Under an assumption of continuity,  the following theorem extends  that result from vectors to  processes.

\begin{Theorem}\label{T3} :  Let $(k(x,y), x,y \in E)$ be a  jointly continuous function on $E\times E$ such that $k(x,x) > 0$ for every $x\in E$.
There exists a positive infinitely divisible process $(\psi_x,x \in E)$  such that for every 
$(\alpha_1, \alpha_2, ...,\alpha_n)$ in $\R^n_+$ and every $(x_1,x_2, ...,x_n)$ in $E^n$,
$$
\e[\hbox{exp}\{-{1\over 2} \sum_{i = 1}^n \alpha_i \psi_{x_i}\}] = | I + \alpha K |^{-1/2}
$$
where  $K = (k(x_i,x_j))_{1\leq i,j \leq n}$

if and only if 
\begin{equation}
\label{1}
k(x,y) = d(x)g(x,y)d(y)
\end{equation}
where $d$ is a positive function and $g$ the Green function of a Markov process.
\end{Theorem}

\b Similarly to the case of vectors,  Remark \ref{T4}.1 leads to the following corollary

\begin{Corollary} : Let $(k(x,y), x,y \in E)$ be a  jointly continuous function on $E\times E$ such that $k(x,x) > 0$ for every $x\in E$. Let $(\psi_x,x \in E)$ be a process such that  for every 
$(\alpha_1, \alpha_2, ...,\alpha_n)$ in $\R^n_+$ and every $(x_1,x_2, ...,x_n)$ in $E^n$,
$$
\e[\hbox{exp}\{-{1\over 2} \sum_{i = 1}^n \alpha_i \psi_{x_i}\}] = | I + \alpha K |^{-1/2}
$$
where  $K = (k(x_i,x_j))_{1\leq i,j \leq n}$.

\b Then $(\psi_x,x \in E)$ is infinitely divisible if and only if it is associated to a Markov process.
\end{Corollary}

\bigskip

\bigskip

 \section{Permanental  random point processes.}\label{S4}

\bigskip

\bigskip

\b The context of this section is  a locally compact Hausdorff space $E$ with a countable basis, 
$\lambda$ is a nonnegative Radon measure on $E$, and $ Q$ is the
space of nonnegative integer-valued Radon measures on $E$. Shirai and Takahashi \cite{ST} have extended the notion of Boson random point process by introducing the following distributions denoted by $\mu_{\alpha,K}$. The corresponding random point processes are sometimes called {\bf permanental point processes} or Bosonian point processes.

\begin{Definition}:  For $K$  a bounded integral operator on $L^2(E,\lambda)$ and $\alpha$ is a
fixed positive number, the distribution $\mu_{\alpha,K}$ on $Q$ satisfies, when it exists
\begin{equation}
\int_{Q}\mu_{\alpha,K}(d\xi) \hbox{exp}(-<\xi,f>) = \hbox{Det}(I + \alpha
K_{\phi})^{-1/\alpha} 
\end{equation}
  for each nonnegative measurable function $f$ on $E$ with compact support where
$K_{\phi}$ stands for the trace class operator defined by
$$ K_{\phi}(x,y) = \sqrt{\phi(x)} K(x,y)  \sqrt{\phi(y)}$$
and $$\phi(x) = 1 - \hbox{exp}(-f(x)).$$
The function Det denotes the Fredholm determinant.
\end{Definition}
For a different presentation of these distributions one can read the paper of Hough et al \cite{HKPV}.

\b When $\alpha = 1$, $\mu_{1,K}$ is the distribution of a the configurations of a Bosonian gaz.
Shirai and Takahashi have established sufficient conditions on the operator $K$ for the existence of the distribution $\mu_{\alpha, K}$. In particular  for $K$ bounded integral operator on $L^2(E,\lambda)$ and $\alpha$ any
fixed positive number, they have established that if $(\alpha,K)$ satisfies 
\smallskip

\b  {\bf (B)} : {\sl the kernel function of the operator $J_{\alpha} = K(I + \alpha K)^{-1}$ is nonnegative}
\bigskip

\b then $\mu_{\alpha,K}$ exists and is infinitely divisible.

\bigskip

\b Looking carefully at their proof, one can actually formulate another condition on $(\alpha,K)$ for the existence of 
$\mu_{\alpha,K}$ that will be more appropriate to our purpose

\begin{Proposition}\label{P2}: For $K$ bounded integral operator on $L^2(R,\lambda)$, $\alpha$ any
fixed positive number, and $S$ any compact subset of $E$, set $J_{\alpha}[S]   = (I + \alpha K_S)^{-1} K_S$ 
where $K_S$ is the restriction of $K$ to $S$. If for every $x_1,x_2,...,x_n$ in $S$, $det_{\alpha} (J_{\alpha}[S](x_i,x_j))_{1\leq i,j\leq n}) \geq 0$
then $\mu_{\alpha,K}$ exists .
\end{Proposition}
Note that $\mu_{\alpha,K}$ is infinitely divisible iff $\mu_{n\alpha,K/n}$ exists for every $n\in \N^*$. Hence if for a given $K$,
 $det_{\beta} (J_{\alpha}[S](x_i,x_j))_{1\leq i,j\leq n}) \geq 0$ for every $\beta > 0$, every $S$ and every $n$, then $\mu_{\alpha,K}$ is infinitely divisible.

\bigskip

\b In the particular case $E = \R^d$, $\alpha = 1$ and $J_1(x,y) = {1\over (4\pi \beta)^{d/2}} \hbox{exp}(- |x-y|^2/4\beta)$ ($K$ is such that  $J_1 = K(I + K)^{-1}$), the distribution $\mu_{1,K}$ can be obtained as  the limit  of the distributions of the positions in $\R^d$ of $N$ identical particles  following the Bose-Einstein statistics in a finite box.  More precisely, one starts from the following random point measure $\mu_{(L,N)}$ which describes the location of an ideal Bosonian gaz composed of $N$ particles in a volume $V = [-L/2, L/2]^d$ with $d \geq  1 $, at a given temperature $T$ 
$$\int_{V^N} \mu_{(L,N)}(d\xi) e^{-(\xi,f)} = cste \int_{V^N} \hbox{exp}(-\sum_{j=1}^N f(x_j)) per(G_L(x_i,x_j))_{1\leq i,j\leq N}dx_1 ...dx_N$$
the constant depends of V and N and is equal to $\int_{V^N} per(G_L(x_i,x_j))_{1\leq i,j\leq N}dx_1 ...dx_N$
where $G_L$ denotes the operator $\hbox{exp}(\beta \Delta_L)$ with $\beta = 1/T$ and $\Delta_L$ is the Laplacian under the periodic condition in $L^2(V)$.

\b As $N$ and $V$ are tending to $\infty$ with $N/V \rightarrow \rho$, \  $\mu_{(L,N)}$ \ converges to a limit depending on $\rho$. Indeed,  denoting by $\rho_c$ the critical density $\int_{\R^d} {dx \over (2\pi)^d} {e^{-\beta |x|^2}\over 1 - 
e^{-\beta |x|^2}}$ which is finite for $d>2$, we have

- if $\rho < \rho_c$, then $\mu_{(L,N)}$ converges to $\mu_{1,K_{\rho}}$, where $K_{\rho} = \l(\rho)J_1( I - \l(\rho) J_1)^{-1}$ and 
${\l}({\rho})$ is a positive constant depending on $\rho$.

\b  This last result provides the  justification to the fact that $\mu_{1,K_{\rho}}$ is the distribution of the configurations of an ideal Bosonian gaz, but the next result is even more illuminating. Indeed  in the case $d>2$ and

- if $\rho \geq  \rho_c$, then  $\mu_{(L,N)}$ converges to a random point process with a distribution $\zeta$  given by 
\begin{equation}
\label{32}
\int_Q \zeta(d\xi) e^{-(\xi,f)} = Det(I + K_{\phi})^{-1} \hbox{exp}\{- (\rho - \rho_c)(\sqrt{1 - e^{-f}}, (I + K_{\phi})^{-1}\sqrt{1 - e^{-f}})\}
\end{equation}
where $K =J_1(I - J_1)^{-1}$
\bigskip

\b The physical explanation of this split convergence is due to the fact that when the density of the gaz becomes higher that $\rho_c$,    a certain proportion of the particles  tend to lower the density by reaching the lowest level of energy. This phenomena, called the  Bose-Einstein condensation, predicted by Einstein in 1925, is intensively studied today specially since this phenomena has been experimentally obtained (for $d=3$ of course) in 1995 by a team at JILA.  It is interesting  to see that the Bose-Einstein condensation phenomena provides an illustration in the case $d=3$ of a  mathematical physics result available for any dimension  $d$ greater than 3.

\b These results have been established by many authors. In particular they are consequences of the works of Bratteli and Robinson \cite{BR} (see Theorem 5.2.32 Chap.5  p.69) and of Fichtner and Freudenberg \cite{FF}. The way Tamura and It\^o 
have obtained these results in \cite{TI1} and \cite{TI}, deserves a special attention because they need neither quantum field theories nor the theory of states on the operator algebras, but mostly a integral formula due to Vere-Jones \cite{VJ}.  Besides
 Tamura and It\^o actually did more than (\ref{32}): in \cite{TI} their proof is based on  the following theorem.

\bigskip

\b {\bf Theorem A} : {\sl  Let $K$ be a bounded symmetric integral operator on $L^2(E,\lambda)$ such that $(1,K)$ satisfies condition {\bf (B)} and  
\begin{equation}
\label{34}
\int_E J_1(x,y) \lambda(dy) \leq 1 \  \ \lambda(dx) \ a.e.
\end{equation}
Then for every $r>0$, there exists a unique random measure with distribution $\zeta_r$ on $Q$  such that for every non-negative measurable function $f$ on $E$
\begin{equation}
\label{33}
\int_Q \zeta_r (d\xi) e^{-<\xi,f>} = \hbox{exp}\{- r(\sqrt{1 - e^{-f}}, (I + K_{\phi})^{-1}\sqrt{1 - e^{-f}})\}
\end{equation}
where $( . , . )$ denotes the inner product of $L^2(E,\lambda)$.
}

\bigskip

\b Tamura and It\^o's result generates several natural remarks. Indeed, in (\ref{32}) the distribution  of the configurations  of the particles is, thanks to Theorem A, the convolution of two distributions: $\mu_{1, K}\ * \zeta_{ \rho - \rho_c}$. 
It is tempting to imagine that  $\mu_{1, K}$ corresponds to   the fraction of the particles with level of energy greater than $0$ 
and  that $\zeta_{ \rho - \rho_c}$ corresponds to the particles that did "coalesce" (ie without kinetic energy or  similarly in a quantic state equal to $1$.) Indeed ${\l}({\rho})$ is a continuous function of $\rho$ on $(0,  \rho_c]$ that takes the value $1$ at $\rho_c$. Hence the distribution of the configurations of particles with density $\rho_c$ has the distribution $\mu_{1, K}$.  The question is: are the configurations of the particles  with $0$ kinetic energy independent of the configurations of the moving particles?  We will   answer  that question.
\smallskip

\b Another natural remark is the following:  The assumption of condition (B) in Theorem A, insures that $\mu_{1,K}$ is infinitely divisible. Besides, since $\zeta_r = (\zeta_{r/n})^{*n}$, $\zeta_r$ is infinitely divisible too. Consequently the distribution $\zeta$ given by (\ref{32}) is also infinitely divisible.  Moreover, thanks to Theorem A, the distribution $\zeta$ exists for any $K$ such that $(1,K)$ satisfies condition (B) and (\ref{34}). That way we hence obtain a family of infinitely divisible distributions.  Who are they? We will give an answer to that question in Theorem \ref{T6}. 

\smallskip

\b Besides, in their paper \cite{ST} (Theorem 6.12) Shirai and Takahashi have obtained a factorization involving $\mu_{\alpha,K}$ for $(\alpha ,K)$ satisfying condition (B) (see (\ref{40}) below). In the case $\alpha = 1$, is this factorization connected to (\ref{32})? We will show that the answer is affirmative and that they are both  direct consequences of the infinite divisibility of $\mu_{\alpha,K}$. 

\smallskip

\b To analyze further the results of Tamura and It\^o, we will use the notion of Cox process.

\begin{Definition}\label{D1}: A Cox process is a Poisson point process with a random intensity  $\sigma$ on the Radon measures on E, hence its distribution $\Pi_{\sigma}$ satisfies  for every nonnegative measurable function $f$ on $E$ with compact support 
$$\int_{Q}\Pi_{\sigma}(d\xi) \hbox{exp}(-<\xi,f>)  = \e[ \hbox{exp}(- \int_E (1 - e^{-f(x)})\sigma(dx))]$$
\end{Definition}
We will work mostly with Cox processes  with random intensity $\psi(x) \lambda(dx)$ where $(\psi(x), x\in E)$ is a positive process such that $\e(\psi(x))$ is a locally bounded function of $x$. Such a Cox process is said driven by $(\psi,\lambda)$. We will shortly denote its distribution by  $\Pi_{\psi, \lambda}$ or $\Pi_{\psi}$ when there is no ambiguity on the measure $\lambda$.
\smallskip

\b If $\psi$ is a permamental process with a kernel  $(K(x,y), x,y \in E)$ and index $\alpha > 0$, we have for every positive function $f$ with compact support
\begin{eqnarray*}
\int_Q \Pi_{(\psi,\lambda)}(d\xi) \hbox{exp}(-<\xi,f>) =  \e[ \hbox{exp}(- \int_E (1 - e^{-f(x)})\psi(x)\lambda(dx))].
\end{eqnarray*}
By dominated convergence, one shows then that  \ \ $\Pi_{(\psi,\lambda)} = \mu_{\alpha,K}$ .
  A priori the couple $(\alpha,K)$ does not satisfy  condition (B) of Shirai and Takahashi, but the existence of $\psi$ requires that for every $\sigma > 0$ the matrix $Q_{\sigma} = ((I + \sigma K)^{-1} K (x_i,x_j))_{1\leq i,j\leq n}$ satisfies : $det_{\alpha}(Q_{\sigma}) \geq 0$. Hence, choosing $\sigma = \alpha$, the sufficient condition of Proposition \ref{P2} is satisfied. 

  \b Note that for a given couple $(\alpha, K)$ the existence of $\mu_{\alpha,K}$ does not guarantee the existence of a permanental process with kernel $K$ and index $\alpha$. Hence not every  permamental point process is a Cox process. 
\bigskip

\b {\bf Remark \ref{D1}.1}: The infinite divisibility of a Cox process with distribution $\Pi_{\psi}$  is not equivalent to  the infinite divisibility of the process $\psi$. Of course the infinite divisibility of $\psi$ implies the infinite divisibility of $\Pi_{\psi}$, but the converse is not true. 
This fact has been stated in 1975 by  Kallenberg \cite{K}(Ex. 8.6, p.58 chp 8- see also Shanbhag and Westcott (1977)\cite{SW}). These references are actually giving examples of real functions $f$ which are not Laplace transforms although the function $\hbox{exp}(f - 1)$ is a Laplace transform. 
 In Section 3, we saw that a permanental process is not always infinitely divisible. Condition (B) of Shirai and Takahashi  allows to put in evidence examples of squared Gaussian processes $\eta^2$ which are not infinitely divisible although 
the Cox process with distribution $\Pi_{\eta^2}$ is infinitely divisible. To this purpose, consider the example of the ideal Bose gaz, where $J_1(x,y) = {1\over (4\pi \beta)^{d/2}} \hbox{exp}(- |x-y|^2/4\beta)$. Thanks to condition (B), $\mu_{1,K}$ is infinitely divisible, and note that 

\b $\mu_{1,K} = \Pi_{{1\over 2}\eta^2,\lambda} * \Pi_{{1\over 2}\eta^2,\lambda}$,  where $(\eta_x, x \in \R^d)$ is a centered Gaussian process with covariance $(K(x,y), x,y \in \R^d)$ and $\lambda$ is the Lebesgue measure on $\R^d$. Hence $ \Pi_{{1\over 2}\eta^2,\lambda}$ is infinitely divisible. Now, making use of Lemma \ref{L1}, we can easily choose $x,y,z \in \R^d$ such that $(\eta^2_x, \eta^2_y, \eta^2_z)$ is not infinitely divisible. Indeed, if $y_i z_i > 0  $ for every $1 \leq i \leq d$ then 
$|y-z|^2 > |y - x|^2 + |z - x|^2$ for $x$ in $\R^d$ with $|x|$ small enough. The matrix $(K(a,b), a,b \in \{x,y,z\})$ has only positive coefficients and  $(K^{-1}(a,b), a,b \in \{x,y,z\})$ has at least one off-diagonal positive coefficient. Consequently this last matrix  can not be an M-matrix. In Lemma \ref{L3}, we will give a characterization of the infinitely divisible $\Pi_{\psi}$ for $\psi$ nonnegative process.

\bigskip

\b  We now introduce some notation useful for next theorem. Let $(\eta_ x, x\in E)$ be a centered Gaussian process with a covariance $(K(x,y), x,y \in E)$. 
We denote by $a$ a point which is not in $E$. One can set $\eta_a = 0$ and then consider the process $(\eta_x , x \in E\cup \{a\})$ by just defining $K(a,a) = K(x,a) = K(a,x) = 0$ for every $x$ in $E$.  Let $(\psi_x, x\in E)$ be a centered Gaussian process with covariance $(K(x,y) + 1, x,y \in E)$. One can similarly consider the process $(\psi_x, x\in E\cup \{a\})$, by noting that $K(a,a) + 1 = K(x,a) + 1 = K(a,x) + 1 = 1$ for every $x$ in $E$. 
For every $\epsilon$,  we define $$\lambda_{\epsilon} = \lambda + \epsilon \delta_a \ .$$
For $(\phi_x, x\in E\cup \{a\})$ positive random process on $E\cup \{a\}$, we denote by  $\Pi_{\phi, \lambda_{\epsilon}}$ the distribution of a Cox process with random intensity $\phi_x\lambda_{\epsilon}(dx)$ on $E\cup \{a\}$. Of course we have: $\Pi_{\eta^2, \lambda_{\epsilon}} = 
\Pi_{\eta^2, \lambda}$.  Without ambiguity $Q$ will denote the space of nonnegative integer-valued Radon measures on $E\cup \{a\}$. With these notations, we can enunciate the following theorem, whose proof is given at the end of this section.

\begin{Theorem}\label{T6} : Let $(\eta_ x, x\in E)$ be a centered Gaussian process with a covariance $(K(x,y), x,y \in E)$. 
Let $(\psi_x, x\in E)$ be a centered Gaussian process with covariance $(K(x,y) + 1, x,y \in E)$. Assume that the distribution 
$\Pi_{{1\over2}\eta^2, \lambda}$ is infinitely divisible, then the following five points are equivalent.
\begin{itemize}

\item[{\bf (i)}] The distribution $\Pi_{{1\over2}(\eta + c)^2, \lambda}$ is infinitely divisible for every constant $c$ in $\R$.

\item[{\bf (ii)}] The distribution $\Pi_{{1\over2}\psi^2, \lambda_{\epsilon}}$ is infinitely divisible for every $\epsilon >0$.

\item[{\bf (iii)}]  For every $r>0$ there exists a random measure with distribution $\nu_r$ on $Q$ such that 
\begin{equation}
\label{23}
\Pi_{{1\over2}(\eta + r )^2, \lambda} \ = \ \Pi_{{1\over2}\eta^2, \lambda} \ *\, \, \nu_r
\end{equation}
Moreover the distribution $\nu_r$ satisfies 
\begin{equation}
\label{24}
\int_Q \nu_r(d\xi) e^{-(\xi,f)} = \hbox{exp}\{- {1\over2}r^2(\sqrt{1 - e^{-f}}, (I + K_{\phi})^{-1}\sqrt{1 - e^{-f}})\}
\end{equation}
where the inner product is taken with respect to the measure $\lambda$.

\item[{\bf (iv)}] For every $\epsilon > 0$, for every $r>0$ there exists a random measure with distribution $\nu_r$ on $Q$ satisfying 
(\ref{23}) and (\ref{24}) but for the measure $\lambda_{\epsilon}$ instead of the measure $\lambda$.

\item[{\bf (v)}] The distribution $\Pi_{{1\over2}(\eta + c)^2, \lambda_{\epsilon}}$ is infinitely divisible for every constant $c$ in $\R$ and every $\epsilon > 0$.

\end{itemize}

\end{Theorem}
{\bf Remark \ref{T6}.1} : For sake of clarity we have stated Theorem \ref{T6} for   squared Gaussian process $\eta^2$, ie a permamental process with a symmetric kernel $(K(x,y), x,y \in E)$ and index $\beta = 2$. But, thanks to Proposition \ref{P4}, a similar theorem holds for a permanental process $\phi$ with a kernel $(K(x,y), x,y \in E)$ and an index $\beta >0$. One has just to 
assume the existence of a permanental process $\psi$ with a kernel $(K(x,y) + 1, x,y \in E)$ and index $\beta$, and replace in the above theorem, the process $(\eta + r)^2$ (resp. $\eta^2$  ) by $(\psi | \psi_a = r^2)$ (resp. $(\psi | \psi_a = 0)$ ). 

\b Note that the distribution $\nu_r$ of (iii) is not necessarily  a Cox process. In view of the results of Sections \ref{S2} and \ref{S3}, this shows clearly the difference between the infinite divisibility of $\Pi_{\phi}$
 and the infinite divisibility of $\phi$. We will see in Section \ref{S7} that this remark extends from permanental processes to nonnegative  processes. 
\bigskip

\b We are now in position to analyze further the results of Tamura and Ito. Remember that in the case $\rho > \rho_c$, the case when a Bose-Einstein condensation occurs, the obtained limit $\zeta$  is equal to : $\mu_{1,K} * \zeta_{\rho - \rho_c}$, where $\zeta_r$ is defined by Theorem A.  First we note that thanks to Theorem \ref{T6}, the existence of  $\zeta_{\rho - \rho_c}$ for every $\rho > \rho_c$ is equivalent to the infinite divisibility of $\mu_{1,K+1}$. To check directly this last property we can, for example, verify that $(1,K+1)$ satisfies condition (B). Indeed we have the  following general proposition which does not require symmetry from $K$. 

\begin{Proposition}\label{P3} :  Let $K$ be a integral operator on $L^2(E,\lambda)$ such that $(1,K)$ satisfies condition (B) and
\begin{equation}
\label{340}
\int_E J_1(y, x) \lambda(dy) \leq 1 \  \ \lambda(dx) \ a.e. 
\end{equation}
Then $(1,K+1)$ satisfies condition (B).
\end{Proposition}
{\bf Proof of Proposition \ref{P3}} : We set $\overline{J}_1 = (K + 1)(I +K + 1)^{-1} $. We have to show that $\overline{J}_1$ has a nonnegative kernel. We denote by  $\1$ the operator on $L^1(E,\lambda)$ with the kernel identically equal to $1$. We have then:  

\b $ I + K + \1 = ( I - J_1)^{-1}  + \1 = (I + \1(I - J_1)) (I - J_1)^{-1}$, 
which leads to 

\b $ \overline{J}_1 = (K + \1) (I - J_1) \left( I + \1(I - J_1)\right)$. Let $f$ be an nonnegative element of $L^2(E , \lambda)$. We set $g = (I + \1(I - J_1))^{-1} f $ and similarly $ f = g +  \1(I - J_1)g$.  Note that $ \1(I - J_1)g$ is a constant that we denote by $c(g)$ and that $c(g) \geq 0$. Indeed, $c(g) =  \1(I - J_1) ( f - c(g))$ and thanks to (\ref{340}) the operator \ $ \1(I - J_1)$ has a positive kernel. Hence  if $c(g) < 0$ then $f-c(g)$ is a positive function and so is 
 $\1(I - J_1) ( f - c(g))$. Consequently $c(g)$ can not be negative. 
 
 \b We have:
$$\overline{J}_1 f =  (K + \1) (I - J_1)g = J_1 g  +   \1(I - J_1)g = J_1(f  - c(g))  + c(g) = J_1f + (I - J_1)c(g).$$
Thanks to  (\ref{340})  $(I - J_1)  c(g) \geq 0$ and thanks to (B), $J_1f $ is nonnegative. Consequently $\overline{J}_1 f $ is non negative. $\Box$ 
\bigskip

\b  Thanks to Theorem \ref{T6}, Proposition \ref{P3} provides  Theorem A of Tamura and Ito.

\b Restraining our attention to the case of ideal Bosonian particles, 
 we see that, with the notations of Theorem \ref{T6}, (\ref{32}) becomes
  \begin{equation}
\label{244}
  \zeta   = \mu_{1,K} * \nu_{\sqrt{2(\rho - \rho_c)}}.
\end{equation}       
Can we interpret $ \nu_{\sqrt{2(\rho - \rho_c)}}$ as the law of the configurations of the particles with $0$ kinetic energy and density $\rho - \rho_c$ ? These particles are at temperature $T = 1/\beta$ and the distribution $\zeta$ depends on $T$. Now imagine that we can lower the temperature $T$ until $0$, we  have then $\rho_c \rightarrow 0$ and for any positive function $f$ with compact support, we easily obtain 
$$\int_Q \zeta(d\xi) e^{-(\xi,f)} \longrightarrow \hbox{exp}\{- \rho  \int_{\R^d}(1 - e^{-f(x)})dx\}$$
The obtained limit is the distribution of a Poisson point process with uniform intensity $\rho dx$ on $\R^d$. But at temperature $T = 0$, all the particles are at  $0$ level of kinetic energy. Hence this  limit  $\Pi_{\rho, dx}$ is the distribution of the configurations of particles, with density $\rho$,  at $0$ level of energy and temperature $0$. This has been already established (differently) by Goldin et al. \cite{GG}.  Now remember that once the $0$  state of kinetic energy is reached, the particles dont move anymore, hence the law of their configurations should not vary when the temperature goes down. But obviously $\nu_{\sqrt{2(\rho - \rho_c)}}$ is different from $\Pi_{(\rho - \rho_c),dx}$. Consequently the answer to the above question is negative.  This implies that the presence of particles in the Bose-Einstein condensation state influences the position of the still moving particles.        
 \smallskip
 
 \b Besides, the expression (\ref{244}) can be rewritten as      
        $$\zeta =  \Pi_{{1\over2}\eta^2} \ *\  \nu_{\sqrt{\rho - \rho_c}}\ *  \ \Pi_{{1\over2}\eta^2} \ * \ \nu_{\sqrt{\rho - \rho_c}}$$
 which leads to 
\begin{equation}
\label{245}
 \zeta = \Pi_{\psi} 
 \end{equation}
  with $(\psi_x, x\in E) \= ({1\over2}(\eta_x + \sqrt {\rho - \rho_c})^2 +  {1\over2}(\tilde{\eta}_x + \sqrt{\rho - \rho_c})^2 , x\in E)$
 where $\eta$ and $\tilde{\eta}$ are  two independent centered Gaussian processes with covariance $(K(x,y), x,y \in \R^d)$. Under this writing  it appears    that $\zeta$ is the distribution of a Cox process.  Similarly (\ref{244}) leads to 
  \begin{equation}
\label{246}
\zeta = \Pi_{{1\over2}(\eta + \sqrt {2(\rho - \rho_c)})^2 } \ *\ \Pi_{{1\over2}\eta^2}
\end{equation}  
Under this last form, physicists can provide an  interpretation in terms of fields (instead of particles). We thank Yvan Castin from Laboratoire Koestler-Brossel for the following  explanation.
The Bosonic field $(\phi(x), x \in \R^d)$ satisfies: $\phi(x) = \phi_0  + \phi_e(x)$, where $\phi_0$ is a (spatially) uniform field corresponding to the condensated particles and $(\phi_e(x), x\in \R^d)$ is the field corresponding to the excited particles. This last field $\phi_e$ is a complex Gaussian field:  $\phi_e = {1\over \sqrt{2}}(\eta + i \tilde{\eta})$. Besides $\phi_0$ is taken uniformly equal to $\sqrt {\rho - \rho_c}$.  The real component of $\phi_e$ can interfere with $\phi_0$ and provides the part $ \Pi_{{1\over2}(\eta + \sqrt {2(\rho - \rho_c)})^2 }$ of $\zeta$. While the complex component of $\phi_e$ does not interfere with $\phi_0$  and its contribution to $\zeta$ is the same as for the gaz without condensation: $\Pi_{{1\over2}\eta^2}$.

 \bigskip
 \b To prove Theorem \ref{T6} we will use the following characterization of the infinite divisible random measure. According  to Theorem 11.2 (chp 11,p.79) in Kallenberg's book \cite{K}, a random measure with distribution $\zeta$ is infinitely divisible  iff there exists for  almost every $x$,w.r.t.  $\e(\zeta)$, a random measure  with distribution $\mu_x$ on $Q$ such that 
\begin{equation}
\label{31}
\zeta^{(x)} \ = \zeta *\,\, \mu_x
\end{equation}
where $(\zeta^{(x)}, x \in E)$ denotes  the Palm measures of $\zeta$.
\bigskip

\b In the special case of a couple $(\alpha,K)$ satisfying condition (B), we hence obtain the existence of $\mu_x$ such that 
\begin{equation}
\label{40}
\mu_{\alpha,K}^{(x)} \ = \mu_{\alpha,K}*\,\, \mu_x
\end{equation}
which is precisely the factorization obtained by Shirai and Takahashi. But here note that it is 
 seen as an immediate  consequence of  the infinite divisibility of $\mu_{\alpha,K}$.
 
 \b We are going to make use of Kallenberg's Theorem (\ref{31}) to characterize the infinitely divisible $\Pi_{\psi}$.

\begin{Lemma}\label{L3} : Let $\Pi_{\psi}$ be the distribution of a Cox process directed by a positive process $(\psi_x, x\in E)$ with respect to $\lambda$.  Assume that $\e(\psi_x)$ is a locally bounded function of $x$.  For each $b$ in $E$ such that $\e(\psi_b)>0$, denote by $(\psi^{(b)}(x) , x\in E)$ the process  $\psi$ under $\e[ {\psi(b)\over \e( \psi(b))};  \ . ]$. The Palm measure at $b$  of $\Pi_{\psi}$, denoted by $\Pi^{b}_{\psi}$,  admits the following factorization for almost every $b$,w.r.t. $\e(\psi_x) \lambda(dx)$
$$\Pi^{b}_{\psi} = \Pi_{\psi^{(b)}} \ * \ \delta_b.$$
The distribution  $\Pi_{\psi}$ is infinitely divisible iff for almost every $b$ w.r.t. $\e(\psi_x) \lambda(dx)$, there exists a random measure with distribution $\mu_b$ such that 
$$ \Pi_{\psi^{(b)}}  \ * \ \delta_b = \Pi_{\psi} \, *\, \mu_b.$$
\end{Lemma}
 {\bf Proof of Lemma \ref{L3} } : For every nonnegative function $f$ on $E$, we have : 
\begin{equation}
\label{16}
\int_Q \Pi_{\psi}(d\xi) e^{-<\xi,f>} = \e[\hbox{exp}(-\int_E (1 - e^{-f(x)}) \psi(x) \lambda(dx))].
\end{equation}
Call $X$ the Cox process with distribution  $\Pi_{\psi}$, then $X$ admits  a first moment measure $M$  on ${\cal B} (E)$ defined by 
$\e(X(A)) = M(A)$, for every $A \in{ \cal B} (E)$ ie 
$M(A) =\e( \int_A \psi(x) \lambda(dx))$. Let $(\Pi_{\psi}^x, x\in E)$ denote the family of Palm measures of $\Pi_{\psi}$, we define  $\tilde{\Pi}_{\psi}^x$ by 
$\Pi_{\psi}^x = \tilde{\Pi}_{\psi}^x   \ * \ \delta_x$. This means that   $(\tilde{\Pi}_{\psi}^x, x\in E)$ satisfies:
$$\int_Q \Pi_{\psi}(d\xi) \int_E \xi(dx) u(\xi,x) = \int_E  M(dx)   \int_Q\tilde{\Pi}_{\psi}^x(d\xi) u(\xi + \delta_x , x)$$
As a consequence of this desintegration formula, we have for any $f$ and $g$ nonnegative functions on $E$ with compact support included in a compact set $A$:
$$-{d\over dt} \int_Q \Pi_{\psi}(d\xi) e^{-<\xi,f + tg>}|_{t = 0} = \int_E g(x) \int_Q \tilde{\Pi}^x_{\psi}(d\xi) e^{-<\xi + \delta_x,f>}M(dx)$$
which thanks to (\ref{16}) leads to
\begin{eqnarray*}
\int_E && g(x) e^{-f(x)} \lambda(dx)\e[\psi(x) \hbox{exp}(-\int_R (1 - e^{-f(y)}) \psi(y) \lambda(dy))] \\
&&= 
 \int_E g(x) \int_Q \tilde{\Pi}^x_{\psi}(d\xi) e^{-<\xi,f>} e^{-f(x)} \e(\psi(x)) \lambda(dx)
\end{eqnarray*}
Consequently: $$\int_Q \tilde{\Pi}^x_{\psi}(d\xi) e^{-<\xi,f>} =  \e[{\psi(x)\over \e(\psi(x))} \hbox{exp}(-\int_R (1 - e^{-f(y)}) \psi(y) \lambda(dy))] $$
hence: $$ \tilde{\Pi}^x_{\psi} =  \Pi_{{\psi}^{(x)}} \ \ \ \e[\psi(x)]\lambda(dx) a.e. . $$
Thanks to (\ref{31}) the infinite divisibility of $\Pi_{\psi}$ is equivalent to the existence  of ${\mu}_b$, $\e[\psi(b)]\lambda(db)$ a.e. such that 
$\Pi_{\psi}^b = \Pi_{\psi} \, *\, {\mu}_b.$
 $\Box$
\bigskip

\b {\bf Proof of Theorem \ref{T6}} :  Let $N$ be a standard Gaussian variable independent of $\eta$. Then $\eta + N$ is a centered Gaussian process with covariance $K + 1$. Hence we take : $\psi = \eta + N$. 
\smallskip

\b {\bf (ii) $\Rightarrow$ (iii)}: 
Assume that  $\Pi_{{1\over2}(\eta + N)^2, \lambda_{\epsilon}}$ is infinitely divisible for every $\epsilon > 0$. 

\b Denote by $\Pi^{(x)}, x \in E\cup \{a\}$ the Palm measures of $\Pi_{{1\over2}(\eta + N)^2, \lambda_{\epsilon}}$. According to Lemma \ref{L3}  there exists for  almost every $x$,w.r.t.  $\e(\eta_x + N)^2 \lambda_{\epsilon}(dx)$, a random measure  with distribution $\mu_x$ on $Q$ such that 
$$\Pi^{(x)} \ = \ \Pi_{{1\over2}(\eta + N)^2, \lambda_{\epsilon}} \, *\,\, \mu_x.$$
Since $ \lambda_{\epsilon}(\{a\}) = \epsilon > 0 $, 
we have $$\Pi^{(a)} \ = \ \Pi_{{1\over2}(\eta + N)^2, \lambda_{\epsilon}}\, *\,\, \mu_a.$$
Besides,  we have $$\Pi^{(a)} \ =  \ \tilde{\Pi}^{(a)} * \delta_a $$
where  $\tilde{\Pi}^{(a)}$ is the law of a Cox process with intensity ${1\over2}(\eta + N)^2$ under $\e(N^2, \ . \ )$ with respect to  $\lambda_{\epsilon}$. Consequently, we obtain 
\begin{equation}
\label{36}
 \tilde{\Pi}^{(a)} \ *\  \delta_a  \ = \Pi_{{1\over2}(\eta + N)^2, \lambda_{\epsilon}}\, *\,\, {\mu}_a
 \end{equation}

\b For a fixed positive constant $r>0$, the finite-dimensional Laplace transforms of the process ${1\over2}(\eta + r)^2$  are given by
$$\e(\hbox{exp}\{-{1\over 2} \sum_{i= 1}^n \alpha_i (\eta_{x_i} + r)^2\} =  
|I + \alpha K|^{-1/2}\hbox{exp}\{- {1\over 2}r^2 1^t(I + \alpha K)^{-1}\alpha 1 \}$$
for every $x_1,x_2,...,x_n$ in $E\cup \{a\}$  where $1 = (1,1,...1)$ belongs to $\R^n$ and $1^t$ denotes its transpose. Consequently for every nonnegative function $f$ on $E\cup \{a\}$
\begin{equation}
\label{371}
\int_Q{\Pi}_{{1\over2}(\eta + r)^2,\lambda_{\epsilon}}(d\xi)e^{-<\xi,f>} =  \e[\hbox{exp}\{-{1\over2}\int_{E\cup \{a\}}(1 - e^{-f(y)}) \eta^2_y \lambda_{\epsilon}(dy)]
\hbox{exp}\{-{1\over2}r^2 F(f,\eta)\}
\end{equation}
with $F(f, \eta) = (\sqrt{1 - e^{-f}}, (I + K_{\phi})^{-1}\sqrt{1 - e^{-f}})$, where the inner product is with respect to $\lambda_{\epsilon}$.  Note that
\begin{equation}
\label{370}
F(f, \eta) = \int_E\sqrt{1 - e^{-f(x)}} (I + K_{\phi})^{-1}\sqrt{1 - e^{-f}}(x) \lambda(dx) + \epsilon \sqrt{1 - e^{-f(a)}}.
\end{equation}
We obtain
\begin{equation}
\label{37}
\int_Q{\Pi}_{{1\over2}(\eta + N)^2,\lambda_{\epsilon}}(d\xi)e^{-<\xi,f>} =  \e[\hbox{exp}\{-{1\over2}\int_{E\cup \{a\}}(1 - e^{-f(y)}) \eta^2_y \lambda_{\epsilon}(dy)]
\e( \hbox{exp}\{-{1\over2}N^2 F(f,\eta)\})
\end{equation}
and similarly
\begin{equation}
\label{38}
\int_{Q}\tilde{\Pi}^{(a)}(d\xi)e^{-<\xi,f>} =  \e[ \hbox{exp}\{-{1\over2}\int_{E\cup \{a\}}(1 - e^{-f(y)}) \eta^2_y \lambda_{\epsilon}(dy)]\e[N^2 \hbox{exp}\{-{1\over2}N^2 F(f,\eta)\}]
\end{equation}
Now, making use of (\ref{37}) and (\ref{38}), equation (\ref{36}) gives 
$$\int_{Q}{\mu}_a(d\xi)e^{<\xi,f>} = \e[N^2 \hbox{exp}\{-{1\over2}N^2 F(f,\eta)\}] \e( \hbox{exp}\{-{1\over2}N^2 F(f,\eta)\})^{-1} \ \hbox{e}^{-f(a)}$$
which thanks to elementary computations on the standard Gaussian law leads to
 $$\int_{Q} {\mu}_a(d\xi)e^{<\xi,f>} = (1 + F(f,\eta))^{-1} \ \hbox{e}^{-f(a)} = \e[e^{-F(f,\eta) T}]\ \hbox{e}^{-f(a)}$$
 where $T$ is an exponential variable  with parameter 1 independent of $\eta$. By multiplying then each member of the above equation by $\int_{Q} \Pi_{{1\over2}\eta^2}(d\xi)e^{<\xi,f>}$, we obtain
$$\int_{Q}{\mu}_a(d\xi)e^{<\xi,f>}  \int_{Q} \Pi_{{1\over2}\eta^2}(d\xi)e^{<\xi,f>} =  \int_{Q} \Pi_{{1\over2}\eta^2}(d\xi)e^{<\xi,f>} \e[e^{-F(f,\eta) T}]\ \hbox{e}^{-f(a)}$$
which implies
$${\mu}_a \ *\  \Pi_{{1\over2}\eta^2} \ =\  \Pi_{{1\over2}(\eta + \sqrt{T})^2, \lambda_{\epsilon}} \ * \ \delta_a \ .$$
Now note that $\Pi_{{1\over2}\eta^2}$ does not charge configurations including the site $a$, hence there exists a distribution $\tilde{\mu}_a$ such that:  ${\mu}_a \ = \tilde{\mu}_a * \delta_a$. We finally obtain 
$$\tilde{\mu}_a \ *\  \Pi_{{1\over2}\eta^2} \ =\  \Pi_{{1\over2}(\eta + \sqrt{T})^2, \lambda_{\epsilon}}.$$
Denote by $X_a$,  $Y_{\eta^2}$ and  $Y_{(\eta + \sqrt{T})^2}$ the random measures corresponding respectively to the distributions $\tilde{\mu}_a$,   $\Pi_{{1\over2}\eta^2}$ and  $\Pi_{{1\over2}(\eta + \sqrt{T})^2, \lambda_{\epsilon}}$. We have
$$X_a \ + \ Y_{\eta^2}\ \= \ Y_{(\eta + \sqrt{T})^2}$$
In particular, since $Y_{\eta^2}(\{a\}) = 0$, we have: $X_a(\{a\})  \= \ Y_{(\eta + \sqrt{T})^2}(\{a\})$. Hence
$$(X_a \ + \ Y_{\eta^2}, X_a(\{a\}))\=(Y_{(\eta + \sqrt{T})^2}, Y_{(\eta + \sqrt{T})^2}(\{a\}).)$$
 Besides, thanks to (\ref{370}), we know that $Y_{(\eta + \sqrt{T})^2}(\{a\}) = N_{\epsilon T}$ where $(N_t, t \geq 0)$ is a Poisson process independent of $({Y_{(\eta + \sqrt{T})^2}}_{| E}, T)$. Similarly $X_a(\{a\}) = N'_{\epsilon T_a}$ where $T_a$ an exponential variable with parameter $1$, and $(N'_t, t \geq 0)$ is a Poisson process independent of $( (X_a \ + \ Y_{\eta^2})_{|E},T_a)$. Hence
 $$((X_a \ + \ Y_{\eta^2})_{| E}, N'_{\epsilon T_a})\=({Y_{(\eta + \sqrt{T})^2}}_{|E}, N_{\epsilon T}).$$
 Since this is true for every $\epsilon > 0$, we obtain
$$({X_a}_{|E} \ + \ Y_{\eta^2}, T_a)\=({Y_{(\eta + \sqrt{T})^2}}_{|E}, T)$$
which leads to 
$$({X_a}_{|E}\ |T_a = r)\ + \ Y_{\eta^2} \= {Y_{(\eta + \sqrt{r})^2}}_{|E} $$
for almost every $r>0$. In terms of distribution, this means that there exists a random measure on $E$ with distribution $\nu_r$ satisfying
$$ \nu_r \ * \ \Pi_{{1\over2}\eta^2} = \Pi_{{1\over2}(\eta +r)^2,\lambda}$$
Thanks to the above equation and to (\ref{371}), we obtain that 
$$\int_Q \nu_r(d\xi) e^{-(\xi,f)} = \hbox{exp}\{- {1\over2}r^2\sqrt{1 - e^{-f}}, (I + K_{\phi})^{-1}\sqrt{1 - e^{-f}})\}$$
the inner product being with respect to the measure $\lambda $. We use now the result contained in Exercise 5.1 p.33 Chap.3 in Kallenberg's book \cite{K}, to check that for any sequence $(r_n, n \in \N)$ of rational numbers converging to a given $r$, the sequence $(\nu_{r_n})$ converges to  a limit distribution satisfying both (\ref{23}) and (\ref{24}) for the measure $\lambda$. Hence (iii) is established for every $r>0$. Since for every real $r$, \ $(\eta+r)^2 \= (\eta - r)^2$, (iii) is obtained for every real $r$.
 \smallskip
 
 \b  {\bf (iii) $\Rightarrow$ (i)}: By assumption $\Pi_{{1\over2}\eta^2,\lambda}$ is infinitely divisible. Besides, since $\nu_{r} = (\nu_{r/\sqrt{n}})^{*n}$, $\nu_r$ is infinitely divisible. Hence $\Pi_{{1\over2}(\eta+ r)^2,\lambda}$ is infinitely divisible for every constant $r$ in $\R$. 
 \smallskip
 
 \b  {\bf (i) $\Rightarrow$ (v)}: Assume that $\Pi_{{1\over2}(\eta + r)^2,\lambda}$ is infinitely divisible for every constant $r$. We have 
 \begin{eqnarray*}
  \int_Q\Pi_{{1\over2}(\eta+r)^2,\lambda_{\epsilon}}(d\xi) e^{-(\xi,f)} 
 &=&  \e[\hbox{exp}\{{1\over2}\int_{E\cup \{a\}}(1 - e^{-f(y)}) (\eta_y+ r)^2 \lambda_{\epsilon}(dy)\}]\\
 &=&  \e[\hbox{exp}\{{1\over2}\int_{E}(1 - e^{-f(y)}) (\eta_y+ r)^2 \lambda(dy)\}] \ \hbox{exp}\{-{1\over2}(1 - e^{-f(a)})\epsilon r^2\}
\end{eqnarray*}
 hence 
 $$\Pi_{{1\over2}(\eta+r)^2,\lambda_{\epsilon}} \ = \  \Pi_{{1\over2}(\eta+r)^2,\lambda} \ * \ \Pi_{{1\over2}r^2, \epsilon \delta_a} \ .$$
 As the convolution of two infinitely divisible distributions, $\Pi_{(\eta+r)^2,\lambda_{\epsilon}}$ is infinitely divisible too for every $\epsilon >0$.
 \smallskip
  
  \b  {\bf (v) $\Rightarrow$ (iv)}:  We keep the notation of the proof of  "(ii) $\Rightarrow$ (iii)".
  \begin{eqnarray*}
  \int_Q\Pi_{{1\over2}(\eta+r)^2,\lambda_{\epsilon}}(d\xi) e^{-(\xi,f)} 
   =  \e[\hbox{exp}\{-{1\over2}\int_{E}(1 - e^{-f(y)}) \eta^2_y \lambda(dy)] \  \hbox{exp}\{-{1\over2}r^2 F(f,\eta)\}
\end{eqnarray*}
For every integer $n$ an every $r$,  $ [\int_Q\Pi_{{1\over2}(\eta+r)^2,\lambda_{\epsilon}}(d\xi) e^{-(\xi,f)}]^{1/n}$ is still a Laplace transform of  a random measure on $Q$. This is true  in particular for $r\sqrt{n}$, but 
 $$ [\int_Q\Pi_{{1\over2}(\eta+r\sqrt{n})^2,\lambda_{\epsilon}}(d\xi) e^{-(\xi,f)}]^{1/n} = 
  \e[\hbox{exp}\{-{1\over2}\int_{E}(1 - e^{-f(y)}) \eta^2_y \lambda(dy)]^{1/n} \  \hbox{exp}\{-{1\over2}r^2 F(f,\eta)\}
 $$
 By letting $n$ tend to $\infty$, we obtain, thanks again to Kallenberg's result (exercise 5.1 p.33 Chap.3 in \cite{K}) that there exists a limit distribution with Laplace transform $\hbox{exp}\{-{1\over2}r^2 F(f,\eta)\}$.
  \smallskip
 
 \b  {\bf (iv) $\Rightarrow$ (ii)}: We start from $\Pi_{{1\over2}(\eta + r )^2, \lambda_{\epsilon}} \ = \ \Pi_{{1\over2}\eta^2, \lambda_{\epsilon}} \ *\, \, \nu_r$.
 By integrating the above equation with respect to $\p(N\in dr)$ we obtain
 $$\Pi_{{1\over2}(\eta+ N)^2,\lambda_{\epsilon}} \ = \ \Pi_{{1\over2}\eta^2, \lambda_{\epsilon}} \ *\, \, \nu_N$$
 where $\nu_N$ denotes the distribution satisfying for every positive function $f$ on $E\cup\{a\}$
 $$\int_Q \nu_N(d\xi) e^{-(\xi,f)} = \int_{\R}\hbox{exp}\{- {1\over2}r^2 F(f,\eta)\}\p(N\in dr)$$
 $F(f,\eta)$ is defined by (\ref{371}).
Since $N^2$ is an infinitely divisible variable, for every integer $n$, there exists an i.i.d. sequence $(Z_1,Z_2,...,Z_n)$ of positive variables such that $ N^2 \= Z_1 + Z_2 + ... + Z_n$. Hence we have
\begin{eqnarray*}
 \int_Q \nu_N(d\xi) e^{-(\xi,f)}& =& \e[\hbox{exp}\{- {1\over2}N^2 F(f,\eta)\}] \\
&=& \e[\hbox{exp}\{- {1\over2}Z_1 F(f,\eta)\}]^n \\
& = & [\int_Q \nu_{\sqrt{Z_1}}(d\xi) e^{-(\xi,f)}]^n
\end{eqnarray*}
consequently $\nu_N$ is infinitely divisible and so is $\Pi_{{1\over2}(\eta+ N)^2,\lambda_{\epsilon}}$ ( ie $\Pi_{{1\over2}\psi^2,\lambda_{\epsilon}}$).   $\Box$

\bigskip

\bigskip

 \section{ Characterization of the infinitely divisible positive processes}\label{S7}
\bigskip

\begin{Lemma}\label{L2} : Let $(\psi_x, x\in E)$ be a positive process.  For every $a$ such that $\e(\psi_a) > 0$, denote by $(\psi_x^{(a)}, x \in E)$ the process having the law of $(\psi_x, x\in E)$ under the probability 
${1\over  \e(\psi_a)} \e(\psi_a,\,  .)$.  Then, $\psi$ is infinitely divisible if and only if for every  $a$ such that $\e(\psi_a) > 0$, there exists a process $(l^{(a)}_x , x \in E)$ independent of $\psi$ such that 
\begin{equation}
\label{17}
\psi^{(a)}\  \= \ \psi +\  l^{(a)}
\end{equation}
\end{Lemma}
\b Thanks to the above lemma,  one could have immediately see in view of Dynkin's isomorphism Theorem, that the Gaussian processes involved there must have an infinitely divisible squared.  Lemma \ref{L2} allows also to measure the difference between the infinite divisibility of  the Cox process driven by $(\psi,\lambda)$ and the infinite divisibility of  $\psi$. 
Indeed if $\Pi_{\psi,\lambda}$ is infinitely divisible and for each $b$ in $E$ the random measure with distribution $\mu_b$ of Lemma \ref{L3} is the addition of a Cox process and a Dirac at $b$, then $\psi$ is also infinitely divisible.

\b Coming back to the isomorphism Theorem and its extension to permanental processes (Theorem \ref{T2}), we see that these relations in law are nothing else but a NSC for the infinite divisibility of the considered positive process $\psi$. Conversely for an infinitely divisible positive process $\psi$, the relation (\ref{17}) can be seen as an isomorphism Theorem between $\psi$ and some process $( (l^{(a)}_x, x\in E), a\in E)$  (note the connection (\ref{30}) - established below - existing between  the different  $l^{(a)}$ as $a$ varies). In that sense we can see Theorem \ref{T6} (and more generally Lemma \ref{L3})  as a super-isomorphism Theorem involving an infinitely divisible Cox process. 

\bigskip

\bigskip

\b {\bf Proof of Lemma \ref{L2}} :  If $\psi$ is infinitely divisible then for every $x = (x_1,x_2,...,x_n) \in E^n$, there exists $\nu_x$ Levy measure on 
$\R^n$ such that $\int_{\R_+^n} 1 \wedge |y| \nu_x(dy) < \infty$ and 
for every $\alpha = (\alpha_1,\alpha_2,...,\alpha_n)$ in $\R_+^n$
\begin{equation}
\label{21}
\e(e^{-\sum_{i = 1}^n \alpha_i \psi_{x_i}}) = \hbox{exp}\{- \int_{\R_+^n} (1 - e^{-(\alpha,y)})\nu_x(dy)\}
\end{equation}
where $(\alpha,y) = \sum_{i = 1}^n \alpha_i y_i$.

\b We hence have: 
$$\e({\psi(x_1)\over \e(\psi(x_1))}e^{-\sum_{i = 1}^n \alpha_i \psi_{x_i}})  
= \e(e^{-\sum_{i = 1}^n \alpha_i \psi_{x_i}})
 \int_{\R_+^n} {y_1\over\e(\psi(x_1))} e^{-(\alpha,y)}\nu_x(dy)$$
 which put in evidence the existence of a process $l^{(x_1)}$ such that $\psi^{(x_1)} = \psi + l^{(x_1)}$.
\smallskip
\b Conversely, assume that for every $a$ there exists a process $l^{(a)}$ satisfying (\ref{17}). By computing the law of $\psi$ under$ \e(\psi_a\psi_b,\ .]$, we see that for every couple $(a,b)$ of $E$, we must have  
\begin{equation}
\label{30}
c_a \e[l_b^{(a)}F(l_x^{(a)}, x \in E)] = c_b \e[l_a^{(b)}F(l_x^{(b)}, x \in E)]
\end{equation}
where we have set: $c_x = \e(\psi_x)$ for every $x$ in $E$. Besides we have: 
$${{\partial\over \partial \alpha_1} \e(e^{-\sum_{i = 1}^n \alpha_i \psi_{x_i}})\over \e(e^{-\sum_{i = 1}^n \alpha_i \psi_{x_i}})} = - \e(e^{-\sum_{i = 1}^n \alpha_i l^{(a)}_{x_i}})\e(\psi_{x_1})$$
hence
\begin{equation}
\label{22}
\e(e^{-\sum_{i = 1}^n \alpha_i \psi_{x_i}}) = \e(e^{-\sum_{i = 1}^n \alpha_i \psi_{x_i}})|_{\alpha_1 = 0}\,
\hbox{exp}\{- \e(\psi_{x_1})\e[{1 - e^{-\alpha_1 l^{(a)}_{x_1}}\over l^{(a)}_{x_1}} e^{-\sum_{i = 2}^n \alpha_i l^{(a)}_{x_i}}]\}.
\end{equation}
For $n=1$, we immediately obtain thanks to (\ref{22}): $\e(e^{-\alpha_1 \psi_{x_1}}) = e^{-\int_0^{\infty}(1 - e^{-\alpha_1y_1}) \nu_{x}(dy_1)}$

\b where $\nu_x(dy_1) = {\e(\psi_{x_1})\over y_1} \p(l^{(a)}_{x_1} \in dy_1)$.

\b Assume now that for $n-1$ the law of $(\psi_{x_1},\psi_{x_2},...,\psi_{x_{n-1}})$ is given by
$$\e(e^{-\sum_{i = 1}^{n-1} \alpha_i \psi_{x_i}}) = \hbox{exp}\{-\int_{[0,\infty)^{n-1}} (1 - e^{-\sum_{i = 1}^{n-1} \alpha_i y_i})\nu_{x}(dy)\},$$
with $\nu_x(dy) = {c_{x_1} \over y_1}  \p(l^{(x_1)}_{x} \in dy)$. We have then thanks to (\ref{30}): 
$\nu_x(dy) = \int_{\R+}{ c_{x_n} \over y_{n} } \p(l^{(x_n)}_{x} \in dy, l^{(x_n)}_{x_n}\in dy_n)$, for every $x_n$ distinct of $x_1,x_2,...,x_{n-1}$. 

\b Now using (\ref{22}), we obtain 
\begin{eqnarray*}
\e(e^{-\sum_{i = 1}^n \alpha_i \psi_{x_i}})& =& \hbox{exp}\{-\int_{[0,\infty)^{n-1}} (1 - e^{-\sum_{i = 2}^{n} \alpha_i y_i})\nu_{x}(dy)\}\\
 &&\hbox{exp}\{-\int_{[0,\infty)^{n}} (e^{-\sum_{i = 2}^n \alpha_i y_i} - e^{-\sum_{i = 1}^{n} \alpha_i y_i}){ c_{x_1} \over y_{1} } \p(l^{(x_1)}_{x} \in dy_1dy_2...dy_n)\}\\
&=  & \hbox{exp}\{-\int_{[0,\infty)^{n}} (1 - e^{-\sum_{i = 1}^{n} \alpha_i y_i}){ c_{x_1} \over y_{1} } \p(l^{(x_1)}_{x} \in dy_1dy_2...dy_n)\}
\end{eqnarray*}

\bigskip

\bigskip

 \section{Proofs of Sections \ref{S1}, \ref{S2} and \ref{S3}}\label{S6}

\bigskip
\b {\bf Proof of Proposition \ref{P4}}: We note first that for any $\alpha = Diag(\alpha_i, 1 \leq i \leq n)$
$$(I + \alpha (G + \delta)) = (I + \alpha G)(I + \delta(I + \alpha G)^{-1} \alpha \1)$$
where $\1$ denotes the $ n\times n$ matrix with all its coefficients equal to $1$. Hence 
$$|I + \alpha (G + \delta)| = |I + \alpha G| | I + D \1|$$
where $D$ is the diagonal matrix such that $D_{ii} =  \sum_{j= 1}^n ( \delta(I + \alpha G)^{-1} \alpha)_{ij}$.
Now note that $|I + D \1| = 1 + Tr(D) = 1 + 1^tD$, denoting by $1^t$ the vector $(1,1,...1)$ of $\R^n$ (see for example
\cite{VJ} identity (4) in section 2). Consequently 
\begin{equation}
|I + \alpha (G + \delta)| = |I + \alpha G|  ( 1 + \delta 1^t (I + \alpha G)^{-1}\alpha 1).
\end{equation}
 Then we obtain
 $$ \e[\hbox{exp}\{-{1\over 2} \sum_{i = 1}^n 
 \alpha_i \psi_{x_i}\}] = | I + \alpha G |^{-1/\beta} (1 + \delta 1^t (I + \alpha G)^{-1}\alpha 1)^{-1/\beta} .$$
A squared centered Gaussian variable is infinitely divisible, hence there exists   a positive random variable $X$, independent of of $\phi$, such that : $\e(e^{-{t\over 2}X}) = (1 + \delta t)^{-1/\beta}$. We obtain 
$$ \e[\hbox{exp}\{-{1\over 2} \sum_{i = 1}^n 
 \alpha_i \psi_{x_i}\}] =  | I + \alpha G |^{-1/\beta} \e[\hbox{exp}\{- {1\over 2}(1^t (I + \alpha G)^{-1}\alpha 1)X\}]$$
Since this identity is true for every $n$, we can write it for $n+1$ with $x_{n+1} = a$. We set $\tilde{G} = (G(x_i,x_j)_{1 \leq i,j \leq n+1}$ and $\tilde{\alpha} = diag(\alpha_i)_{1\leq i \leq n+1}$. 
We choose to take $x_{n+1} = a$.  Note that we have then: $1^t (I + \tilde{\alpha} \tilde{G})^{-1}\tilde{\alpha} 1 = \alpha_{n+1}\  + \ 1^t (I + \alpha G)^{-1}\alpha 1$,
where $1^t$ and $1$ denote without ambiguity   vectors of $\R^{n+1}$ in the left hand side and of $\R^{n}$ in the right hand side.  Hence we obtain
\begin{equation}
\label{06}
 \e[\hbox{exp}\{-{1\over 2} \sum_{i = 1}^{n+1} 
 \alpha_i \psi_{x_i}\}] = | I + \alpha G |^{-1/\beta}
 \e[\hbox{exp}\{- {1\over 2}(1^t (I + \alpha G)^{-1}\alpha 1 +  \alpha_{n+1})X\}] 
 \end{equation}
which first gives :  $\e[\hbox{exp}\{-{1\over 2} \alpha_{n+1} \psi_{a}\}]  =  \e[\hbox{exp}\{- {1\over 2} \alpha_{n+1}X\}]$.
Developing then (\ref{06}), we obtain
\begin{eqnarray*}
&&\int_0^{+\infty} \p(\psi_a \in dr) e^{-{1\over 2} \alpha_{n+1}r}\e[\hbox{exp}\{-{1\over 2} \sum_{i = 1}^{n} 
 \alpha_i \psi_{x_i}\} | \psi_a = r] \\
 && =  \e[\hbox{exp}\{-{1\over 2} \sum_{i = 1}^n 
 \alpha_i \phi_{x_i}\}] \int_0^{+\infty} \p(\psi_a \in dr) e^{-{1\over 2} \alpha_{n+1}r}
 \hbox{exp}\{- {r\over 2}(1^t (I + \alpha G)^{-1}\alpha 1\}
 \end{eqnarray*}
which leads to (\ref{02}).
\bigskip

\b {\bf Proof of Proposition \ref{P5} }: First note that by assumption $\tilde{G} = {1\over 2}( G + G^t)$ is positive definite and hence $G$ is invertible. We have for any  diagonal matrix $\alpha$ with non-negative coefficients:
 \begin{eqnarray*}
 |I + (\alpha - A) \tilde{G}| 
  &=& | \tilde{G}| |  \tilde{G}^{-1} - A + \alpha| \\
  &=& | \tilde{G}| |  \tilde{G}^{-1} - A| | I + \alpha( \tilde{G}^{-1} - A)^{-1}| \\
\end{eqnarray*}
which leads to
\begin{equation}
\label{061}
  |I + (\alpha - A) \tilde{G}|  =  |\tilde{G}| | G|^{-1} | I + \alpha G|
 \end{equation}
Besides, for the complex Gaussian vector  $\Lambda = (\eta_{x_j} + i \tilde{\eta}_{x_j})_{1 \leq j \leq n}$ where $\tilde{\eta}$ is an independent copy of $\eta$, we have
$$ |I + (\alpha - A) \tilde{G}|^{-1} = \e[\hbox{exp}\{-{1\over 2} < (\alpha - A)\Lambda, \overline{\Lambda}>\}]$$
and in particular for  $\alpha = 0$: $\e[\hbox{exp}\{{1\over 2} <  A\Lambda, \overline{\Lambda}>\}]
= { |G|\over |\tilde{G}|} $.  Thanks to (\ref{061}), we obtain then  (\ref{05}).

\b Note that for the Gaussian vector $\eta = (\eta_{x_j})_{1 \leq j\leq n}$, we have 
$$ \e[\hbox{exp}\{-{1\over 2} < (\alpha - A)\eta, {\eta}>\}] = | I + (\alpha - {1\over 2}(A + A^t) )\tilde{G}|^{-1/2}$$
hence (\ref{05}) which connects the norm of $(\psi, \tilde{\psi})$ to the norm of $(\eta, \tilde{\eta})$, can not be reduced to 
 a one dimensional identity that would connect $\psi^2$ to $\eta^2$.

\bigskip

\b {\bf Proof of Theorem \ref{T1} } :  We use Proposition 4.5 of Vere-Jones paper \cite{V}. First we note that thanks to Assertion $(D_{16})$ in Chap. 6, p.135 of Berman and Plemmons's book \cite {BP} , all the real eigenvalues of $G$ are positive. Then, since 
the resolvent matrices $Q_{\sigma} = \sigma G (I + \sigma G)^{-1}$, $\sigma > 0$, have only nonnegative coefficients, they are all $\beta$-positive definite for every $\beta > 0$. Hence $\psi$ is well defined.
 \bigskip
 
 \b {\bf Proof of Theorem \ref{T2}} : Thanks to Theorem 1, we have  for  $(x_1,x_2,...x_n) \in E^n$:
 \begin{equation}
\label{8}
\langle\hbox{exp}\{-{1\over 2} \sum_{i = 1}^n \alpha_i \psi_{x_i}\}\rangle = | I + \alpha G |^{-1/2}.
\end{equation}
Hence setting $x_1 = a$, we obtain : 
 \begin{equation}
 \label{6}
 \langle\psi_a\hbox{exp}\{-{1\over 2} \sum_{i = 1}^n \alpha_i \psi_{x_i}\}\rangle = {\partial\over \partial \alpha_1}(| I + \alpha G |) \ | I + \alpha G |^{-3/2}.
\end{equation}
A development of  $| I + \alpha G |$ with respect to its first column, gives :
\begin{eqnarray*}
| I + \alpha G |& =& (1 + \alpha _1 g(a,a))  ( I + \alpha G )^{11} - \alpha_1g(x_2,a)( I + \alpha G )^{21}\\
                                & +& \alpha_1g(x_3,a)( I + \alpha G )^{31}+ ... + (-1)^{n+1}g(x_n,a)\alpha_1( I + \alpha G )^{n1}
\end{eqnarray*}
and hence
\begin{eqnarray*}
{\partial\over \partial \alpha_1}(| I + \alpha G |) & =&  g(a,a)  ( I + \alpha G )^{11} - g(x_2,a)( I + \alpha G )^{21}\\
                                & +& g(x_3,a)( I + \alpha G )^{31}+ ... + (-1)^{n+1}g(x_n,a)( I + \alpha G )^{n1}\\
                                & = & |A|
\end{eqnarray*}
 where the matrix $A = (A_{ij})_{1\leq i,j\leq n}$ is defined by $A_{ij} = (I + \alpha G )_{ij}$ if $i\not=1$ and 
$A_{1j} = g(x_j,a)$.
Consequently (\ref{6}) becomes:
\begin{equation}
 \label{7}
 \langle\psi_a\hbox{exp}\{-{1\over 2} \sum_{i = 1}^n \alpha_i \psi_{x_i}\}\rangle = |A|  \ | I + \alpha G |^{-3/2}.
\end{equation}
Besides, we know (see for example Marcus and Rosen \cite{MR1} Lemma 2.6.2) that 
$$
 \tilde{\p}_{a}( \hbox{exp}\{- \sum_{i = 1}^n \alpha_i  L^{x_i}_{\infty}) =  {|A|\over g(a,a)  | I + \alpha G |} 
$$
which  together with (\ref{8}) and  (\ref{7}) give Theorem 2. $\Box$
\bigskip

\b{\bf Proof of Corollary \ref{C1}} :  For every $(x_1,x_2,...x_n) \in E^n$ such that $x_1 = a$,  the L\'evy measure of $\psi$ satisfies 
$$\langle \hbox{exp}\{-{1\over 2} \sum_{i = 1}^n \alpha_i \psi_{x_i}\}\rangle  = \hbox{exp}\{-\int_{\R^n}(1 -   \hbox{e}^{-\sum_{i = 1}^n \alpha_i {y_i}})\nu(dy)\}. $$
Hence, we obtain
$$
\langle \psi(a) \hbox{exp}\{-{1\over 2} \sum_{i = 1}^n \alpha_i \psi_{x_i}\}\rangle  =  \langle \hbox{exp}\{-{1\over 2} \sum_{i = 1}^n \alpha_i \psi_{x_i}\}\rangle\ \int_{\R^n}(2y_1  \hbox{e}^{-\sum_{i = 1}^n \alpha_i {y_i}})\nu(dy) .
 $$
Consequently, thanks to Theorem 2
$$g(a,a) \tilde{\p}_{a}( \hbox{exp}\{- \sum_{i = 1}^n \alpha_i  L^{x_i}_{\infty}) = \int_{\R^n}(2y_1  \hbox{e}^{-{1\over 2}\sum_{i = 1}^n \alpha_i {y_i}})\nu(dy) 
 $$
which leads to
\begin{eqnarray*}
\nu_{(\psi(a)/2,\psi(x_2)/2,...\psi(x_n)/2)}\!\!\!\!\!\!&&\!\!\!\!\!\!(dy_1,dy_2,...,dy_n) \\
&=& {g(a,a)\over 2y_1} \tilde{\p}_{a}(L^a_{\infty}\in dy_1,
L^{x_i}_{\infty}\in dy_i, 2\leq i\leq n).
\end{eqnarray*}
\bigskip

\b {\bf Proof of Corollary \ref{C2}} : We denote by $g_{ \tau_{S_{\theta}}}$ the Green function of $X$ killed at $\tau_{S_{\theta}}$, and by $g_{T_a}$ the Green function of $X$ killed at $T_a$. It has been proved in \cite{EK1}  that for every $x,y\in E$ : $g_{ \tau_{S_{\theta}}} (x,a) = g_{ \tau_{S_{\theta}}} (a,x) = 1/\theta$ and $g_{ \tau_{S_{\theta}}} (x,y)  = g_{T_a}(x,y) + 1/\theta$.
Hence  for the process $X$  killed at $\tau_{S_{\theta}}$:  $\tilde{\p}_{a} = \p_a$, and by Theorem \ref{T2} : 
\begin{equation}
\label{14}
{\p}_{a}\langle F(  L^x_{\tau_{S_{\theta}}}  + {1\over 2} \psi_x ; x\in E) \rangle = \langle  {\psi_a\over g(a,a)} F({1\over 2} \psi_x ; x\in E) \rangle
\end{equation}
with $\langle\hbox{exp}\{-{1\over 2} \sum_{i = 1}^n \alpha_i \psi_{x_i}\}\rangle = | I + \alpha G_{ \tau_{S_{\theta}}} |^{-1/2}$
where $G_{ \tau_{S_{\theta}}} = (g_{T_a}(x_i,x_j) + 1/\theta)_{1\leq  i,j \leq n}$. We set $G_{T_a} =  (g_{T_a}(x_i,x_j) )_{1\leq  i,j \leq n}$.

\b Thanks to Proposition \ref{P4} we can define a measurable function ${f}$ on $\R_+^n\times E^n$ such that  for any $\alpha$ and any $x = (x_1,x_2,...,x_n) \in E^n$
\begin{equation}
\label{35}
| I + \alpha G_{ \tau_{S_{\theta}}} |^{-1/2} = | I + \alpha G_{T_a}|^{-1/2}( 1  + {1\over\theta}{f}(\alpha,x))^{-1/2}.
\end{equation}
Assuming that $x_1 = a$, we obtain then  
\begin{eqnarray*}
-1/2 \langle  {\psi_a} \hbox{exp}\{-{1\over 2} \sum_{i = 1}^n \alpha_i \psi_{x_i}\}\rangle& = &{\partial \over \partial \alpha_1} ( | I + \alpha G_{T_a}|^{-1/2}( 1  + {1\over\theta}{f}(\alpha,x))^{-1/2})\\
&=& -{1\over 2\theta} | I + \alpha G_{T_a}|^{-1/2}{\partial f(\alpha,x)\over \partial \alpha_1}( 1  + {1\over\theta}{f}(\alpha,x))^{-3/2}
\end{eqnarray*}
since  $| I + \alpha G_{T_a}|$ does not depend of $\alpha_1$. 

\b Besides, since the process$(L^._{\tau_r}, r>0)$ is a L\'evy process under $\p_a$, there exists a measurable function $h$ on $\R_+^n\times E^n$ such that for every $(\alpha,x)$
$$\p_a ( \hbox{exp}\{-{1\over 2} \sum_{i = 1}^n \alpha_i L^{x_i}_{\tau_r}\}) = \hbox{e}^{- h(\alpha,x)r}.$$
Hence (\ref{14}) can be expressed as
$$\e(\hbox{e}^{h(\alpha,x)S_{\theta}})| I + \alpha G_{T_a}|^{-1/2}( 1  + {1\over\theta}{f}(\alpha,x))^{-1/2}  =  | I + \alpha G_{T_a}|^{-1/2}{\partial f(\alpha,x)\over \partial \alpha_1}( 1  + {1\over\theta}{f}(\alpha,x))^{-3/2}$$
which is equivalent to 
$$( 1  + {1\over\theta}{h}(\alpha,x))^{-1}     = {\partial f(\alpha,x)\over \partial \alpha_1}( 1  + {1\over\theta}{f}(\alpha,x))^{-1}.$$
Consequently for every $\theta > 0$, we have : $ {\partial f(\alpha,x)\over \partial \alpha_1} = { \theta  + {f}(\alpha,x)\over  \theta  +{h}(\alpha,x)}$. By letting $\theta $ tend to $\infty$, we finally obtain  $ {\partial f(\alpha,x)\over \partial \alpha_1} = 1$ and $f(\alpha, x) = h(\alpha, x)$. We can now rewrite (\ref{35}) as follows
$$\langle\hbox{exp}\{-{1\over 2} \sum_{i = 1}^n \alpha_i \psi_{x_i}\}\rangle =  \p_a \langle \hbox{exp}\{- \sum_{i = 1}^n \alpha_i ( L^{x_i}_{\tau_{S_{\theta}}} + {1\over 2}\phi_{x_i}\}\rangle.$$
By conditionning on both sides by the respective value at $a$ of the processes, we obtain:
$$\langle\hbox{exp}\{-{1\over 2} \sum_{i = 1}^n \alpha_i \psi_{x_i}\}|\psi_a = r \rangle =  
\p_a \langle \hbox{exp}\{- \sum_{i = 1}^n \alpha_i ( L^{x_i}_{\tau_{S_{\theta}}} + {1\over2}\phi_{x_i})\} | L^a_{\tau_{S_{\theta}} }= r\rangle$$
which is equivalent to 
$$\langle\hbox{exp}\{-{1\over 2} \sum_{i = 1}^n \alpha_i \psi_{x_i}\}|\psi_a = r \rangle =  
\p_a \langle \hbox{exp}\{- \sum_{i = 1}^n \alpha_i ( L^{x_i}_{\tau_r} + {1\over2}\phi_{x_i})\}\rangle.$$
$\Box$
\bigskip

 \bigskip
\b {\bf Proof of Lemma \ref{L1}} :

\b First assume (\ref{03}).
For $\alpha = diag(\alpha_1, \alpha_2, ..., \alpha_n)$, we set:  $F(\alpha) = | I + \alpha G |^{-1}$.
For $a >0$, we set: $Q_{a} = a G(I + a G)^{-1}$  and for $S = diag(s_i)_{1\leq i \leq n}$
with $|s_i| \leq 1 , 1 \leq i\leq n$ we define :
$P_{a }(S) = |I-Q_{a}| \ |I - Q_{a}S|^{-1}$.

\b If the function $F(\alpha) $ is the Laplace transform of an  infinitely divisible vector $(\psi_1,\psi_2,...,\psi_n)$, then for every $a >0$,  the function $P_{a }(S)$ is the probability generating function of an infinitely divisible vector. This has been used in the symmetric case by Griffiths \cite{G} and Griffiths and Milne \cite{GM}, but  it is still true without assumption of symmetry. Indeed, we have :  $P_{a }(S)  = F( a (I - S))$, which can be rewritten as
$$P_a(S) = \e[ \Pi_{i = 1}^n s_i ^{N_i}]$$
where conditionnally to $(\psi_1,\psi_2,...,\psi_n)$, $N_1,N_2,...,N_n$ are $n$ independent Poisson variables with respective parameter $a\psi_1,a\psi_2,...,a\psi_n$.

\b Griffiths and Mile \cite{GM} have established the following criterion

\bigskip

\b {\bf   Theorem B} : {\sl  Let $Q$ be a $n\times n$ real matrix. The function $|I-Q| \ |I - QS|^{-1}$ is an infinitely divisible probability generating function if and only if

(i) the eigenvalues of $Q$ are strictly bounded in modulus by $1$;

(ii) $Q_{ii}\geq 0$ and $Q_{ij}Q_{ji} \geq 0, i\not=j,  i,j\in \{1,2,...,n\}$;

(iii) for every $k\leq n$, for every subset $\{i_1,i_2,...,i_k\}$ of $k$ distinct indices from $\{1,2,...,n\}$ 
$$ T_{i_1,i_2}T_{i_2,i_3}ÉÉÉT_{i_{k-1},i_k}T_{i_k,i_1} \geq 0$$
where  $T = Q + Q^t$\ ( $Q^t$ denotes the transpose of $Q$).
}
\bigskip

\b We have : 
\begin{equation}
\label{3}
I - Q_{a} = {1\over a} ( a^{-1} I  + G)^{-1}
\end{equation}
For every $a >0$, $P_{a}(S)$ is an infinitely divisible probability generating function. We choose $a$ large enough in order  that 
 for every $(i,j)$ 
 
\b  if $G^{-1}_{ij} \not= 0$, then 
$G^{-1}_{ij}$ and $( a^{-1} I  + G)^{-1}_{ij}$ have the same sign 

\b if $G^{-1}_{ij} \not= 0$ and $G^{-1}_{ji} = 0$, then
$(G_{ij}^{-1} +  G_{ji}^{-1})$ and $( a^{-1} I  + G)_{ij}^{-1} + ( a^{-1} I  + G^t)_{ji}^{-1})$ have the same sign.

\b Thanks to Theorem B (ii),  and (\ref{3}), we obtain
\begin{equation}
\label{5}
G^{-1}_{ij} G^{-1}_{ji}\geq 0
\end{equation}
Making use of the argument of Bapat to prove Theorem 1 \cite{B},  we know that there exists a signature matrix $\sigma$ such that the off-diagonal terms of the matrix $\sigma (-T)\sigma$ are all negative.
\begin{equation}
\label{4}
I - T = {1\over a} \{ ( a^{-1} I  + G)^{-1} + ( a^{-1} I  + G^t)^{-1} \}
\end{equation}
We have   $(\sigma (I -T)\sigma)_{ij} \leq 0$ for $i\not=j$, hence we obtain
$$\sigma(i)\sigma(j)\{ ( a^{-1} I  + G)_{ij}^{-1} + ( a^{-1} I  + G)_{ji}^{-1} \} \leq 0.$$

\b  if $G^{-1}_{ij} G^{-1}_{ji}\not= 0$, then $\sigma(i)\sigma(j)(G_{ij}^{-1} +  G_{ji}^{-1})  \leq 0$
and hence thanks to (\ref{5}), we obtain $\sigma(i)\sigma(j)G_{ij}^{-1} \leq 0$.

\b if $G^{-1}_{ij} \not= 0$ and $G^{-1}_{ji} = 0$, then
$\sigma(i)\sigma(j)G_{ij}^{-1} \leq 0$.

\b Consequently for every $(i,j)$, $\sigma(i)\sigma(j)G_{ij}^{-1} \leq 0$. Making use of the result $(N_{38})$  of Berman and Plimmons in \cite{BP} p.137 chap.6, 
we finally obtain that $\sigma G^{-1} \sigma$ is an $M$-matrix. 
\smallskip

\b Conversely, assume that there exists a signature matrix $S$ such that $S G^{-1} S$ is an $M$-matrix.  We can then reproduce  the proof of Theorem 3.2 \cite{EK} to show that
$$S(i)G(i,j)S(j) = d(i) g(i,j) d(j)$$
 where $d$ is a deterministic function on $\{1,2,...,n\}$ and $g$ is the Green function of a transient Markov process 
 with a state space equal to $\{1,2,...,n\}$.  
 
 \b Using then Theorem  \ref{T1} and Corollary \ref{C1}, we conclude that  (\ref{03}) is satisfied.
 
 \bigskip
 
 \b {\bf Proof of Theorem \ref{T4} } : Thanks to Theorem   \ref{T1} and Corollary \ref{C1}, the sufficiency of this condition is immediat. The necessity follows from the argument developed in \cite{EK} to establish Theorem 3.2. This argument is based on Bapat's criterion. 
 
\smallskip 
  
 \b {\bf Proof of Theorem \ref{T3}} : The sufficiency of the condition (\ref{1}) follows from Theorem \ref{T1} and Corollary \ref{C1}. To prove the necessity,  we can, thanks to Lemma \ref{L1},
 make use of the proof of Theorem 3.4 \cite{EK} which works similarly since there was no use of symmetry in the proof.
 \bigskip

\bigskip

\bigskip

 \section{A conjecture on  $\alpha$-permanents }\label{S5}

\b Here is Shirai and Takahashi's conjecture:  

\smallskip

\b {\sl  Let $\alpha$ be in $[0,2]$. Then $\hbox{det}_{\alpha}A$ is nonnegative for
every nonnegative definite square matrix $A$.}
\smallskip

\b Shirai and Takahashi have made substancial progresses in the direction of proving the conjecture. We send the interested reader to the last section of their paper \cite{ST}.  Besides in \cite{HKPV}, the authors give an example of a $3\times 3$ positive definite matrix $K$ such that $det_{\alpha} (K) < 0$ for $\alpha > 4$. 

\b Actually one  easily checks that the conjecture is true for $3\times 3$ positive definite matrices. 
We would like just to point out the fact that in view of the results of Vere-Jones, this conjecture has the following  more appealing form for a probabilist: 
\smallskip
 
 \b {\sl For every centered Gaussian vector $(\eta_1,\eta_2,...,\eta_n)$ for every $\delta \geq 1$
 $$\e[\hbox{exp}\{- \sum_{i = 1}^n z_i \eta^2_i\}]^{\delta}$$
 is still a Laplace transform in $(z_1,z_2,...,z_n)$.}

\vspace{1 cm}

\noindent
\begin{tabular}{lll}  & Nathalie Eisenbaum& Haya Kaspi\\
   & Laboratoire de Probabilit\'es& Industrial Engineering and
Management\\
  &Universit\'e Paris VI - CNRS &Technion\\
  &4 place Jussieu &  Technion City\\ &75252 Paris Cedex 05,
France&Haifa 32000, Israel
\\ &nae@ccr.jussieu.fr&iehaya@tx.technion.ac.il\\ & &\\ & & \\ & &

\\
\end{tabular}

\end{document}